\documentclass[12pt]{amsart}
\newcommand{\df}{\stackrel{\rm{def}}{=}}

\newcommand{\la}{\lambda}

\newcommand{\K}{\mathcal{K}}

\newcommand{\T}{\mathbb{T}}

\newcommand{\e}{\varepsilon}

\newtheorem{thm}{Theorem}[section]

\newcommand{\ci}[1]{_{ {}_{\scriptstyle #1}}}

\newcommand{\cD}{\mathcal{D}}
\newcommand{\R}{\mathbb{R}}

\numberwithin{equation}{section}

\begin{document}

\title[Bellman functions and inequalities for Haar multipliers]{Bellman functions and 
two weight inequalities for Haar multipliers}
\author{F. Nazarov, S. Treil, A. Volberg}
\address{Department of Mathematics, Michigan State University,
East Lansing, Michigan, 48824}
\email[Nazarov]{fedja@math.msu.edu}
\email[Treil]{treil@math.msu.edu}
\email[Volberg]{volberg@math.msu.edu}
\curraddr[Volberg]{Mathematical Sciences Research Institute,
1000 Centennial Drive, Berkeley, CA 94707-5070}

\thanks{Partially supported by the NSF
grant DMS 9622936, binational Israeli-USA grant BSF 00030, and research programs at MSRI in the Fall of 1995 and in the Fall of 1997.}

\begin{abstract}
We are going to give necessary and sufficient conditions for two weight norm inequalities
for Haar multipliers operators and for square functions. We also give sufficient
conditions for two weight norm inequalities for the Hilbert transform.
\end{abstract}

\maketitle

\setcounter{section}{-1}
\section{Introduction}
\label{intro}

Weighted norm inequalities for singular integral operators appear 
naturally in many areas of analysis,  probability, operator theory ect. 

The one-weight case is now pretty well understood, and the answers 
are given by the famous Helson--Szeg\"{o} theorem and the 
Hunt--Muckenhoupt--Wheden Theorem. The fist one state that the Hilbert 
Transform $H$ is bounded in the weighted space $L^{2}(w)$ if and only if 
$w$ can be represented as $w=\exp\{u+Hv\}$, where $u,\ v\in 
L^{\infty}$, $\|u\|_{\infty}<\pi/2$. 

The Hunt--Muckenhoupt--Wheden Theorem states that the Hilbert 
transform $H$ is bounded in $L^{p}(w)$ if and only if the weight $w$ 
satisfies the so-called Muckenhoupt $A_{p}$ condition
\begin{equation}
\sup_{I} \Bigl( \frac1{|I|} \int_{I} w \Bigr) \cdot \Bigl(\frac1{|I|} 
\int_{I} w^{-1/(p-1)} \Bigr)^{p-1} < \infty,  \tag{\mbox{$A_{p}$}}
\end{equation}
where the {\em supremum} is taken over all intervals $I$. 
This condition is also necessary and sufficient fo boundedness for a 
wide class of singular integral operators, as well as for the 
boundedness of the maximal operator $M$,
$$
Mf(x) = \sup_{I\ni x} \frac1{|I|} \int_{I} |f|  ; 
$$
here supremum is taken over all intervals $I$ containing $x$. 

It is worth mentioning, that there in no direct proof of equivalence 
the Helson--Szeg\"{o} condition and the Muckenhoupt condition 
$A_{2}$. 

Two weight inequalities, i. e. the problem when an operator acts 
from $L^{2}(w) $ to $L^{2}(v)$ (one can also consider $L^{p}$ case, 
even with different exponents $p$, but the $L^{2}$ case is complicated 
enough, so we restrict our attention on it) appears naturally in many areas like 
the theory of Hankel and Toeplitz operators, perturbation theory, etc.

Things look much more complicated in the two-weight case, and it is 
probably an agreement now that there is no simple (Muckenhoupt type) necessary 
and sufficient condition of boundedness of the Hilbert Transform. 

It was a big surprise when Eric Sawyer \cite{Sa1} found necessary and 
sufficient condition for a maximal operator $M$ to be a bounded 
operator from $L^{2}(w)$ to $L^{2}(v)$: his theorem states that it is 
enough to test the boundedness on a very special class of test 
functions, namely only on functions $\chi\ci I 
w^{-1}$ (and we should do the same for the adjoint operator).\footnote{Sawyer's theorem states   more, and 
treats the $L^{p}$-case as well, but we are not going into details 
here} We will call such type of conditions {\em Sawyer type} conditions. 

There is also a two weight analog of Helson--Szeg\"{o} theorem due to M. 
Cotlar and C. Sadosky, see \cite{CS1}. Their approach (which can
be referred as Generalized Bochner Theorem) provides both integral representation and
extension of forms and kernels invariant under the shift operator. Being applied to a
special bilinear form built with the help of the Hilbert transform and two measures,
this approach gives a necessary and sufficient condition for the Hilbert transform to be
bounded between $L^2$-spaces with respect to these measures (see
[CS1]). The approach of  M.Cotlar, C.Sadosky is very intersting because it provides 
a direct link between the lifting theory of Sz.-Nagy and Foias (and thus the scattering
theory) and the continuity of the Hilbert transform in weighted spaces 
(see \cite{S}).

But there is no  analog of the Muckenhoupt $A_{2}$ condition for two 
weights, which is necessary and sufficiend for the boundedness of the 
Hilbert transform. There are quite a few sufficient conditons, let us 
mention a very nice and simple one due to Dechao Zheng 
\cite{Zh}.

In this paper we are going to consider an operator (more precisely, a 
family of operators, the so called Haar multipliers) which can serve as 
a good ``model'' for singular integral operators.  For such operator 
we give necessary and sufficient conditions (of Sawyer type) of the 
boundedness. 

 Our operators appear to  be simpler than the Hilbert transform and we
believe that our approach splits the difficulties of two weight singular integral
estimates and allows to treat these difficulties separately.

So let us explain what is our ``model'' operator.  Let $\mathcal{D}$ 
denote the set of dyadic subarcs of the real line $\R$.  Let $\sigma = 
(\sigma\ci{I})\ci{I\in\cD}$ be a sequence of signs $\pm$.  We will be 
dealing actually with the following family of operators.  Let $I_{-}, 
I_{+}$ denote the left and the right halves of a dyadic interval $I$ 
and let
$$ h\ci{I} = \left\{ \begin{array}{l}
+ |I|^{-1/2} \chi\ci{I_{-}}, \\[.3cm]
- |I|^{-1/2} \chi\ci{I_{+}}.
\end{array} \right. $$
denote a Haar function normalized in $L^{2}=L^{2}(R)$. Let $(\cdot, \cdot)$ denote the scalar
product in $L^2$. We are interested in the following
question: How to describe the pairs
$(\mu,\nu)$ of weights on $\R$ such that all operators $T_{\sigma}$,

%

\begin{equation}
\label{0.1}
T_{\sigma}f = \sum_{I\in D} \sigma\ci{I}(f, h\ci{I})  h\ci{I} 
\end{equation}
are uniformly bounded from $L^{2}(\mu)$ to $L^{2}(\nu)$ with respect to
all possible choice of $\sigma$?

It is easy to see that the measure $\nu$ has to be absolutely 
continuous. It is also not difficult to see that singular part of 
$\mu$ does not help: if operators $T_{\sigma}$ are uniformly bounded 
from $L^{2}(\mu)$ to $L^{2}(\nu)$, then the same holds if we 
replace $\mu$ by its absolutely continous part. So, without loss of 
generality  one can assume that the measures $\mu$ and $\nu$ are 
absolutely continuous, $d\mu=udt$, $d\nu=vdt$. 

If the operators $T_{\sigma}$ are uniformly bounded, then the operators
$$
f\mapsto (f, h\ci I)\ci{L^{2}} h\ci I, \qquad I\in\cD
$$
are uniformly bounded as well. For a fixed $I$ the norm of the above operator can be easily 
computed --- it is just a rank one operator --- and it is equal to 
$\langle v \rangle\ci I^{1/2} \langle u^{-1} \rangle\ci I^{1/2}$. So 
we get a simple necessary condition
$$
\sup_{I\in\cD} \langle v \rangle\ci I \langle u^{-1} \rangle\ci I < \infty,
$$
which can be considered as a two-weight analog of the Muckenhoupt $A_{2}$ 
condition. Unfortunately this condition is not sufficient. 

It is convenient to denote $w:=u^{-1}$. In this notation the unifrom 
boundedness of operator $T_{\sigma}L^{2}(u) \to L^{2}(v)$ is 
equivalent to the uniform boundedness of the operators 
$M_{v}^{1/2} 
T_{\sigma} M_{w}^{1/2}$ in usual (non-weighted) $L^{2}$; here 
$M_{v}$ and $M_{w}$  denote operators of multiplication 
by $v$ and $w$ respectively.

And as we have shown above, the following condition
$$\sup_{I \in \mathcal{D}}\langle v \rangle\ci{I} \langle w \rangle\ci{I} \leq C
$$
is necessary for the uniform boundedness of all $M_{v}^{1/2} 
T_{\sigma} M_{w}^{1/2}$ in $L^{2}$ or equivalently for the uniform 
boundedness of all operators $T_{\sigma}: L^{2}(w^{-1}) \to L^{2}(v)$.

The above uniform boundedness of $T_{\sigma}$ admits a simple 
geometric interpretation in terms of so-called \emph{multipliers}. 
Consider the family ${\mathcal M}$ of bounded operators $A :
L^{2}(w^{-1}) \to L^{2}(v)$ such that they commute with all 
operators $f\mapsto (f,h\ci I) h\ci I$, $I\in\cD$. 
We call this family the family of Haar multipliers.

These are operators given by a simple formula $Ah\ci{I} = a\ci{I}h\ci{I}$. If
$a = \{a\ci{I}\}\ci{I\in \mathcal{D}}$ is given, let us call this operator
$A_{a}$,

Now one can note that uniform boundedness of $T_{\sigma}$ is equivalent to
the inclusion
\begin{equation}
\label{0.3}
 \ell^{\infty} \subset {\mathcal M} 
\end{equation}
in the sense that $A_{a} \in {\mathcal M}$ for all $a \in \ell^{\infty}$.

We are going to investigate the question when \eqref{0.3} holds, that is when
the family of Haar multipliers contains $\ell^{\infty}$.

We are going to formulate three main results now. Notice that Theorems 
\ref{t0.2}
and \ref{t0.3} together
give the necessary and sufficient conditions for
\begin{equation}
\label{0.4}
 \sup_{\sigma} \|M_{v}^{1/2} T_{\sigma} M_{w}^{1/2}\| < \infty 
 \end{equation}
and Theorem \ref{t0.1} also gives the necessary and sufficient 
conditions in a completely diferent form.

\bigskip

\begin{thm}
\label{t0.1}
The family of singular integrals $M_{v}^{1/2} T_{\sigma} M_{w}^{1/2}$ is uniformly
bounded in $L^{2}$ if and only if
\begin{enumerate}
\item  ${\displaystyle \forall \; J \in \mathcal{D} \quad \sup_{\sigma} \frac{1}{|J|} \int_{J}
|T_{\sigma}(\chi\ci{J}w)|^{2} v \, dx \leq C \langle w\rangle\ci{J}\,;}$ 

\item ${\displaystyle \forall \; J \in \mathcal{D} \quad \sup_{\sigma} \frac{1}{|J|} \int_{J}
|T_{\sigma}(\chi\ci{J}v)|^{2} w \, dx \leq C\langle v\rangle\ci{J}\,.}$
\end{enumerate}
\end{thm}
Note, that any of the above conditions 1 or 2 immediately implies taht
$$
\forall \; J \in \mathcal{D} \quad \langle v\rangle_{J}\langle 
w\rangle_{J} \leq C < \infty, 
$$
which is the necessary condition we discussed above.

To formulate the next theorem, let us introduce some notation. 
Let us denote 
$$
\alpha\ci{I} = \left|\dfrac{\langle
v\rangle\ci{I_{-}}-\langle v\rangle\ci{I_{+}}}{\langle v\rangle\ci{I}}\right|
\left|\dfrac{\langle w\rangle\ci{I_{-}}-\langle w\rangle\ci{I_{+}}}{\langle 
w\rangle\ci{I}}\right|
|I|,
$$ 
where $I_{-}, I_{+}$ are left and right halves of $I$.
Consider an integral operator $T_{0}$ given by the
formula
$$
T_{0}f = \sum_{I\in \mathcal{D}} \frac{1}{|I|}\langle f 
\rangle\ci{I}\chi\ci{I}\alpha\ci{I}
$$
whose kernel $k(x,y) = \sum_{I\in \mathcal{D}} {|I|^{-2}}
\chi\ci{I}(x)\chi\ci{I}(y)\alpha_{I} $ is evidently positive.

\bigskip

\begin{thm}
\label{t0.2} The family of operators $M_{v}^{1/2} T_{\sigma} M_{w}^{1/2}$ is uniformly
bounded in $L^{2}$ if and only if the following four assertions hold
simultaneously:
\begin{enumerate}
\item ${\displaystyle \forall \; J \in \mathcal{D} \quad \langle 
v\rangle\ci{J}\langle w\rangle\ci{J} \leq C <
\infty}$;
\item ${\displaystyle \forall \; J \in \mathcal{D} \quad \frac{1}{|J|} \sum_{I\subseteq J}
|\langle v\rangle\ci{I_{-}}-\langle v\rangle\ci{I_{+}}|^{2}\langle 
w\rangle\ci{I}|I| \leq C
\langle v\rangle\ci{J}}$;

\item ${\displaystyle \forall \; J \in \mathcal{D} \quad \frac{1}{|J|} \sum_{I \subseteq J}
|\langle w\rangle\ci{I_{-}}-\langle w\rangle\ci{I_{+}}|^{2}\langle 
v\rangle\ci{I}|I| \leq C
\langle w\rangle\ci{J}}$;
\item  The operator $T_{0}$
 is bounded from $L^{2}(w^{-1})$ to $L^{2}(v)$, or, equivalently, the 
 operator $M_{v}^{1/2} T_{0} M_{w}^{1/2}$ is 
bounded in $L^{2}$.
\end{enumerate}
\end{thm}

\begin{thm}
\label{t0.3}
The operator $T_{0}$ is bounded from $L^{2}(w^{-1})$ to
$L^{2}(v)$ if and only if
\begin{enumerate}
\item ${\displaystyle \forall \; J \in \mathcal{D} \quad \frac{1}{|J|} 
\int_{J} \Bigl(
\sum_{I\subset
J} \frac{1}{|I|} \chi\ci{I}\langle w\rangle\ci{I}\alpha\ci{I}\Bigr)^{2} v \, dx \leq C\langle
w\rangle\ci{J}}$\,;

\item  ${\displaystyle \forall \; J \in \mathcal{D} \quad \frac{1}{|J|} 
\int_{J}\Bigl(\sum_{I\subset J}
\frac{1}{|I|} \chi\ci{I} \langle v\rangle\ci{I} \alpha\ci{I}\Bigr)^{2} w \, dx \leq C\langle
v\rangle\ci{J}}$\,.
\end{enumerate}
\end{thm}

Theorem 0.1 looks surprising. If it would concern an operator with positive kernel it would
be in the vein of Sawyer's weighted theorems from [S1], [S2]. In fact, exactly as in [S1],
[S2], Theorem 0.1 claims that (the family of) integral operators are bounded  if and
only if (the family of) integral operators are bounded on test functions $\chi_J$ and the
adjoint operators are bounded on test functions $\chi_J$. However, the family $T_{\e}$ models
a singular integral operator rather than a positive kernel integral operator. Unlike the case
of positive kernel integral operators it now seems surprising that boundedness
on $\chi_J$ implies boundedness on smaller positive functions. 

On the other hand, the classical $T1$ theorem says exactly the same: if an operator with
Calder\'{o}n-Zygmund kernel is uniformly bounded on $\chi_J$ then it is bounded. Recently, 
(see [NTV] and [T]) the $T1$ theorem was extended to nonhomogeneous spaces (spaces with
non-doubling measure). It turned out that this generalization plays an important role in the
treatment of analytic capacity problems including a famous problem of Ahlfors-Vitushkin.

The method we use to prove these theorems consists of constructing the Bellman function of the
problems we consider. Roughly speaking, we try to solve an extremal problem associated to a
given problem (simply speaking,we try to consider the worst possible case). This leads to
Bellman function of the problem. This approach resembles the approach of Burkholder ([Bu]).
But there is a difference. We are not solving the extremal problem mentioned above (we would
only wish). Instead, we are looking for a kind of subsolution of an associated system of
Partial Differential (In)Equalities.
 
Theorem \ref{t0.3} can  be most probably proved using the combination of ideas frim the
papers of Kalton, Verbitsky [KV] and Sawyer, Wheeden [SW]. Following [KV] we can
introduce the new metric $d(x,y) := \frac{1}{k(x,y)}$, where $k$ is the kernel of $T_0$
and was written above. It is clearly a metric, and balls in this metric are just all
dyadic intervals. Then kernel $k$ satisfies the regularity condition from [SW]. Using
the main result of [SW] we could have given an alternative proof of the theorem.
Unfortunately it is not clear why the new metric space has certain regularity properties.
For example, in [KV] one requires the property that all annuli be nonempty. This is false
in our new metric space. However, it is not clear to us how essential are those
regularity properties of the metric space for the application of the technique of [SW].

{\bf Acknowledgements}. We are grateful to Peter Jones, Robert Fefferman and Igor
Verbitsky for valuable discussions of this paper.

\section{Necessary conditions}
\label{nec-cond}

As it was shown  above in the Introduction, the uniform boundedness 
of $M^{1/2}_v T_\sigma M^{1/2}_w$ implies that 
the operators $f\mapsto (f, h\ci I) h\ci I$ are uniformly bounded as
operators from $L^2(w^{-1})$ to
$L^2(v)$, and the later condition is equivalent to  
\begin{equation}
\label{1.1}
\langle v\rangle\ci{I} \langle w\rangle\ci{I}
\leq C < \infty.
\end{equation}
It is well known that for the case of one weight ($w^{-1}=v$ in our
notation) this caondition is just the famous Muckenhoupt $A_2$
condition, and it is a sufficient condition for the unoform
boundedness of $M^{1/2}_v T_\sigma M^{1/2}_w$.  

But in general case we have other  simple necessary conditions which are independent
of \eqref{1.1}. To get on of the condition let us apply the operator 
$M^{1/2}_v T_\sigma M^{1/2}_w$ to the test function $w^{1/2} \chi\ci
J$. We get
$$
 \int_{\R} \Bigr|\sum_{I\in\mathcal{D}} \sigma\ci{I}
(w \cdot \chi\ci J,h\ci{I})h\ci{I}(x) \Bigl|^{2}v(x) \, dx
$$

Let us now take the average over all possible choices of signs
$\sigma\ci I$.

Fix a function $g$ in $L^{2}(w^{-1})$ and let $\{\e_{I}(\omega)\}$ be
the sequence of independent random variables assuming values $\pm 1$
with probabilities 1/2, 1/2. Then
$$ \int_{\Omega} d\mathbb{P} \int_{\T} |\sum_{I\in\mathcal{D}} \e_{I}(\omega)
(g,h_{I})h_{I}(x)|^{2}v \, dx \leq C\|g\|^{2}_{L^{2}(w^{-1})} $$
and
$$ \int_{\T} \sum_{x\in I} |(g,h_{I})|^{2} \frac{1}{|I|} v(x) dx
\leq C \int_{\T} |g|^{2}w^{-1} $$
which is
$$ \sum_{I \in \mathcal{D}} |\langle g \rangle_{I_{-}} - \langle g \rangle_{I_{+}}|^{2} \langle v
\rangle_{I}
\cdot |I|
\leq C
\int_{\T} |g|^{2}w^{-1} $$
Choose $g = \chi_{J}w$. Then we come to the following necessary
condition (notice that we keep only the summation over $I \subset J$)
\begin{equation}
\forall \, J \in \mathcal{D}, \;\frac{1}{|J|} \sum_{I \subset J}
|\langle w\rangle_{I_{-}}-\langle w \rangle_{I_{+}}|^{2} v_{I} \cdot |I| \leq C \langle w
\rangle_J.
\end{equation}
Symmetrically (we can interchange $w$ and $v$)
\begin{equation}
\forall \, J \in \mathcal{D}, \; \frac{1}{|J|}\sum_{I \subset J}
|\langle v \rangle_{I_{-}} - \langle v \rangle_{I_{+}}|^{2} w_{I} \cdot |I| \leq C \langle v
\rangle_J
\end{equation}
If $ w^{-1} = u = v $, then, say, (1.2) becomes
\begin{equation}
\frac{1}{|J|}\sum_{I\subset J} \langle w^{-1}\rangle_{I} \langle w \rangle_{I}^{2}
\left(\frac{\langle w \rangle_{I_{-}}- \langle w \rangle_{I_{+}}}{\langle w
\rangle_{I}}\right)^{2} |I|
\leq C\langle w \rangle_J
\end{equation}
which follows from $(A_{2})$ condition (1.1) on $w$ (nontrivially, see
e.g. [FKP] or [B]).

But in general neither (1.2) nor (1.3) follows from (1.1). This is shown in [N]. In what
follows we will try to establish to what extent (1.1), (1.2), (1.3) is the full list of
necessary and sufficient conditions.

In general this is not the case, but first we separate some
cases when (1.1)-(1.3) or their simple modifications are sufficient for
\begin{equation}
\sup_{\e} \|T_{\e}\|_{L^{2}(w^{-1})\to L^{2}(v)} \leq C < \infty .
\end{equation}

\setcounter{equation}{0}

\section{2. Reduction of Theorem 0.1 to Therems 0.2, 0.3 and the proof of Theorem 0.2.}

Let us use the notation $(\cdot,\cdot), \; (\cdot,\cdot)_{w^{-1}}$, and
$(\cdot,\cdot)_{v}$ for the dualities in $L^{2}, \, L^{2}(w^{-1})$, and
$L^{2}(v)$. Let $T$ be an arbitrary operator from $L^{2}(w^{-1})$ to
$L^{2}(v)$ with the norm $\|T\|_{L^{2}(w^{-1}) \to L^{2}(v)} =
\sup_{\|f\|_{2} \leq 1,\|g\|_{2} \leq 1} \;
|(Tw^{1/2}g, v^{1/2}f)|$, where $\|\cdot\|_{2}$ denotes the usual
$L^{2}$-norm.

We are interested in estimating from above $\|T\|$. Thus we can assume
for the time being that $f$ and $g$ are bounded  as
long as the final estimates of $\|T\|$ will not depend on these bounds.

Adopting this convention let us write for arbitrary $f, g \in L^{\infty}$ the
following decompositions: the decomposition of $fw^{1/2}$ and $gv^{1/2}$ with respect to
the usual Haar basis: 

\begin{eqnarray*}
fw^{1/2} & = & \sum (f, w^{1/2} h_{J}) h_{J}, \\[.2cm]
gv^{1/2} & = & \sum (g, v^{1/2} h_{I}) h_{I}.
\end{eqnarray*}

The reader can notice that we are using the usual Haar system
$\{h_{I}\}$ and its biorthogonal system $\{ \mbox{const}_{I}
\frac{h_{I}}{v}\}$ in $L^{2}(v)$ to set up the decomposition in the
output space $L^{2}(v)$ and we proceed similarly in the input space $L^{2}(w^{-1})$.

Now the expression $(Tw^{1/2}g, v^{1/2}f)$ involved in the formula for
the norm $\|T\|$ becomes

\begin{equation}
(Tw^{1/2}f, v^{1/2}g) = \sum(g, v^{1/2}
h_{I})(f,w^{1/2}h_{J})(Th_{J},h_{I}).
\end{equation}

For the operators $T_{\e}$ we have

\begin{equation}
(T_{\e}w^{1/2}f,v^{1/2}g) = \sum_{I\in \mathcal{D}} \e_{I}
(g,v^{1/2}h_{I})(f,w^{1/2}h_{I}).
\end{equation}

To estimate the latter expression we are going to use

\subsection{Disbalanced Haar functions}

This is the system of functions $\{h^{}_{I}\}_{I\in \mathcal{D}}$ having
the following properties:

\begin{enumerate}
\item[1)] $h^{w}_{I}$ vanishes outside of $I$ and equals to two
different constants on the left and on the right halves of $I$,
\item[2)] $\int h^{w}_{I} w \, dx = 0$,
\item[3)] $\|h^{w}_{I}\|_{L^{2}(w)} = 1$.
\end{enumerate}

Then this is an orthonormal system in $L^{2}(w)$. Such kind of system
with nonpositive weight $z'(s)ds$ has been used in [CJS] to give a
simple solution of a problem of Calder\'{o}n. The following identity
plays an important part below:

\begin{equation}
x_{I} \cdot h_{I} = h^{w}_{I} + A_{I} \cdot \chi_{I},
\end{equation}
where the constants $x_{I}$ and $A_{I}$ are uniquelty defined by the
properties of $\{h^{w}_{I}\}$ listed above. Let us compute them. Using 2)
and 3) we get
$$ 
x^{2}_{I} \cdot \langle w \rangle_{I} = \|h^{w}_{I}\|^{2}_{L^{2}(w)} + A^{2}_{I}|I| \langle
w\rangle_I = 1 + A^{2}_{I} |I|\langle w \rangle_I 
$$
Considering the scalar product of both parts of (1.8) with the constant
function in $L^{2}(w)$ we get $x_{I}(1,h_{I})_w = A_{I}|I| \langle w \rangle_I$.

Thus

\begin{eqnarray*}
x_{I} & = & \sqrt{\frac{\langle w \rangle_I}{\langle w \rangle_{I_{-}}\langle w
\rangle_{I_{+}}}}; \\[.2cm]
A_{I} & = & \frac{x_{I} }{2\sqrt{|I|}} \frac{\langle w \rangle_{I_{-}} -\langle w
\rangle_{I_{+}}}{\langle w \rangle_{I}}; \\[.2cm]
\frac{A_{I}}{x_{I}} & = & \frac{1}{2 \sqrt{|I|}}
\frac{\langle w \rangle_{I_{-}} -\langle w
\rangle_{I_{+}}}{\langle w \rangle_{I}}.
\end{eqnarray*}
And having this in mind let us plug (2.3) into (2.2) to get 4 sums:

\begin{eqnarray*}
\lefteqn{(T_{\e}w^{1/2}f,v^{1/2}g) = \sum_{I\in\mathcal{D}} \e_{I}
\left(\frac{g}{v^{1/2}},h^{v}_{I}\right)_{v}(fw^{-1/2},h^{w}_{I})_{w}
\frac{1}{x^{v}_{I}x^{w}_{I}} } \\[.2cm]
&& \mbox{} \hspace*{.2in} + \sum_{I\in\mathcal{D}} \e_{I} \left(
\frac{g}{v^{1/2}},h^{v}_{I}\right)_{v}(fw^{1/2})_{I} \frac{1}{x^{v}_{I}}
A^{w}_{I} \frac{1}{x^{w}_{I}} |I| \\[.2cm]
&& \mbox{} \hspace*{.2in} + \sum_{I\in\mathcal{D}} \e_{I} (gv^{1/2})_{I}
(fw^{-1/2},h^{w}_{I})_{w} \cdot \frac{1}{x^{v}_{I}} A^{v}_{I}
\frac{1}{x^{w}_{I}} |I| \\[.2cm]
&& \mbox{} \hspace*{.2in} + \sum_{I\in\mathcal{D}} \e_{I} (gv^{1/2})_{I} \cdot
(fw^{1/2})_{I}
\frac{1}{x^{v}_{I}} A^{v}_{I} |I| \cdot \frac{1}{x^{w}_{I}}
A^{w}_{I} |I| \\[.2cm]
&& \mbox{} = \Sigma_{1} + \Sigma_{2} + \Sigma_{3} + \Sigma_{4}.
\end{eqnarray*}
 
To reduce Theorem 0.1 to Theorems 0.2, 0.3 we will prove that the assumptions of Theorem 0.1
imply the assumptions of Theorems 0.2, 0.3. We already saw that the averaging over signs of
the second assumption of Theorem 0.1 implies (1.2) and
the averaging over signs of the third assumption of Theorem 0.1 implies (1.3). So in our
reduction we can  already use (1.2), (1.3) freely.

Four sums above define four linear operators whose bilinear forms are given by the
corresponding sums. Let us call them $D_{\e}, \Pi_{\e}, \Pi^{'}_{\e}, T_{0,\e}$
correspondingly.

A very important remark now is the following. Suppose that all assumptions of Theorem 0.1 are
satisfied. 
Then operators
$v^{1/2}T_{\e}w^{1/2}$ are uniformly (in $\e$) bounded on
$f,g$ if and only if $T_{0,\e}$ are uniformly bounded on $f,g$.

This is true because of the following analysis of operators $D_{\e}, \Pi_{\e}, \Pi^{'}_{\e}$,
which says that these operators are uniformly (in $\e$) bounded if the assumptions of Theorem
0.1 are satisfied. To analyse these operators is the same as to analyse their bilinear forms
given by $\sum_{i}, i = 1,2,3$.
 
\subsection{Estimates of $\sum_{i}, i = 1,2,3$}

There is nothing to estimate in $\sum_{1}$. Notice that
$\frac{1}{x^{v}_{I}} \leq \sqrt{\langle v \rangle_{I}}$ so we can
use the bound $\frac{1}{x^{v}_{I}}\frac{1}{x^{w}_{I}}\leq \sqrt{\langle v
\rangle_{I}}\sqrt{\langle w \rangle_{I}} \leq C$ to write
\begin{equation}
\Sigma_{1} \leq C \|\frac{g}{v^{1/2}}\|_{L^{2}(v)} \cdot
\|fw^{-1/2}\|_{L^{2}(w)} = C\|g\|_{2} \cdot \|f\|_{2}.
\end{equation}

To estimate $\sum_{2}$ let us notice that
$$ \left| \frac{1}{x_{I}^{v}} A^{w}_{I}
\frac{1}{x^{w}_{I}}\right| |I|
\leq \sqrt{\langle v \rangle_{I}}  \frac{\langle w \rangle_{I_{-}}-\langle w
\rangle_{I_{+}}}{\langle w \rangle_{I}}
\sqrt{|I|}.
$$
Then
\begin{eqnarray*}
\Sigma_{2} & \leq & \left\|\frac{g}{v^{1/2}}\right\|_{L^{2}(v)} \left(
\sum_{I\in\mathcal{D}} \langle fw^{1/2} \rangle_{I}^{2} \langle v \rangle_{I}
\left(\frac{\langle w \rangle_{I_{-}}-\langle w
\rangle_{I_{+}}}{\langle w \rangle_{I}}\right)^{2}|I|\right)^{1/2}
\\[.2cm]
& = & \|g\|_{2} \cdot \left(\sum_{I\in\mathcal{D}} \langle fw^{1/2}\rangle^{2}_{I}\langle v
\rangle_{I}
\left(\frac{\langle w \rangle_{I_{-}}-\langle w
\rangle_{I_{+}}}{\langle w \rangle_{I}}\right)^{2}|I|\right)^{1/2}.
\end{eqnarray*}

To estimate the last expression let us use the following lemma.

\bigskip

{\bf Lemma 2.1.} {\it Let $\{\alpha_{I}\}_{I \in \mathcal{D}}$ be a sequence of
nonnegative numbers. Then $$\sum_{I\in\mathcal{D}}
\langle fw^{1/2}\rangle^{2}_{I}
\alpha_{I}
\leq C
\|f\|^{2}_{2}$$ if and only if for all $ J,\linebreak \frac{1}{|J|}\sum_{I\subset J}
\langle w\rangle_{I}^{2} \alpha_{I} \leq C \langle w \rangle_J$.} 

\bigskip

We postpone the proof till Section 6. If we use this lemma and (1.2) we can finish the
estimate of $\sum_{2}$:

\begin{equation}
\Sigma_{2} \leq C_{2} \|f\|_{2} \cdot \|g\|_{2}.
\end{equation}

Similarly, (1.3) gives

\begin{equation}
\Sigma_{3} \leq C_{3} \|f\|_{2} \cdot \|g\|_{2}.
\end{equation}

We conclude that if the assumptions of Theorem 0.1 hold then, for any given pair $f,g$ of 
$L^2$-functions $\sup_{\e}|(v^{1/2}T_{\e}w^{1/2}f,g)| \leq A \|f\|_2 \|g\|_2$ if and only if 
$\sup_{\e}|(T_{0,\e}f,g)| \leq B \|f\|_2 \|g\|_2$. But this supremum can be computed as
follows

\begin{eqnarray*}
\sup_{\e}\Sigma_{4} = \sup_{\e}\sum \e_{I} \langle gv^{1/2}\rangle_{I}\langle
fw^{1/2}\rangle_{I}
\frac{\langle v \rangle_{I_{-}}-\langle v \rangle_{I_{+}}}{\langle v
\rangle_{I}}\frac{\langle w \rangle_{I_{-}}-\langle w \rangle_{I_{+}}}{\langle w
\rangle_{I}}|I| & = &
\\[.2cm]
\sum_{I\in\mathcal{D}} |(gv^{1/2})_{I}| \cdot |(fw^{1/2})_{I}|
\frac{|\langle v \rangle_{I_{-}}-\langle v\rangle_{I_{+}}|}{\langle v\rangle_{I}} \cdot
\frac{|\langle w \rangle_{I_{-}}-\langle w\rangle_{I_{+}}|}{\langle w\rangle_{I}} |I| 
\end{eqnarray*}

Clearly the estimate

$$
\sum_{I\in\mathcal{D}} |\langle gv^{1/2}\rangle_{I}| \cdot |\langle fw^{1/2}\rangle_{I}|
\frac{|\langle v \rangle_{I_{-}}-\langle v\rangle_{I_{+}}|}{\langle v\rangle_{I}} \cdot
\frac{|\langle w \rangle_{I_{-}}-\langle w\rangle_{I_{+}}|}{\langle w\rangle_{I}}|I| 
 \leq 
C \|f\|_{2} \|g\|_{2}
$$
is equivalent to
$$|(T_{0}f,g)| \leq  C
\|f\|_{2}
\|g\|_{2}.
$$

Now we can finish the reduction. Second and third assumptions of Theorem 0.1 say exactly that
for pairs $f = \chi_{J}w^{1/2}, g = \chi_{J}g$ and $ f = \chi_{J}f, g = \chi_{J}v^{1/2}$ the
supremum of bilinear forms $\sup_{\e}|(v^{1/2}T_{\e}w^{1/2}f,g)| $ has the desired estimate.
The consideration written above then implies that on such pairs the bilinear form
$|(T_{0}f,g)|$ has the desired estimate. But this means exactly that the second and third
assumptions of Theorem 0.3 are satisfied. The reduction is finished. We saw that the
assumptions of Theorem 0.1 imply the first three assumptions of Theorem 0.2 and all
assumptions of Theorem 0.3, but they, in their turn, imply the fourth assumption of Theorem
0.2. So we are left to prove theorems 0.2 and 0.3.

Notice that Theorem 0.2 is already proved. In fact, we established that the estimate
$\sup_{\e}|(T_{\e}f,g)| \leq  C
\|f\|_{2}
\|g\|_{2}
$ implies (1.1)--(1.3). We also know that (1.1)--(1.3) guarantee the equivalence of the
inequality $\sup_{\e}|(T_{\e}f,g)| \leq  C
\|f\|_{2}
\|g\|_{2}
$ with the inequality $|(T_{0}f,g)| \leq  C
\|f\|_{2}
\|g\|_{2}
$. So Theorem 0.2 is completely proved.

\bigskip

\setcounter{equation}{0}

\section{3. $A_{\infty}$ conditions and sufficient conditions.}

In this section we show that relatively simple conditions (1.1)--(1.3) are already
sufficient for uniform boundedness of 
$\|T_{\e}\|_{L^2(w^{-1}\to L^2(v)}$ 
if weights have certain
$A_{\infty}$ properties. In the first theorem below we assume that either $v$ or $w$ belong to
$A_{\infty}$.

\bigskip

{\bf Theorem 3.1} {\it If necessary conditions (1.2), (1.3) hold and one of
the functions $v$ or $w$ is in $A_{\infty}$ then (2.7) holds too.} 

\bigskip

{\bf Proof.} Let $v \in A_{\infty}$. Then ( see [FKP]) we have

\begin{equation}
\forall \; J, \quad \;\frac{1}{|J|} \sum_{I\subset J} \langle v\rangle_{I}
\left(\frac{\langle v\rangle_{I_{-}}-\langle v\rangle_{I_{+}}}{\langle
v\rangle_{I}}\right)^{2} |I|
\leq C \langle v\rangle_J.
\end{equation}

Application of Lemma 2.1 now proves that

\begin{equation}
\sum_{I\in\mathcal{D}} |\langle gv^{1/2}\rangle_{I}|^{2} \cdot \frac{1}{\langle v\rangle_{I}}
\left(
\frac{\langle v\rangle_{I_{-}}-\langle v\rangle_{I_{+}}}{\langle v\rangle_{I}}\right)^{2}|
|I|
\leq C
\|g\|_{2}.
\end{equation}

Then the left part of (2.7) can be estimated as 

$$
 \left(\sum|\langle gv^{1/2}\rangle_{I}|^{2} \frac{1}{\langle v\rangle_{I}}
\left(
\frac{\langle v\rangle_{I_{-}}-\langle v\rangle_{I_{+}}}{\langle v\rangle_{I}}\right)^{2}
I\right)^{1/2}
$$
$$
\left(
\sum|\langle fw^{1/2}\rangle_{I}|^{2} \left(
\frac{\langle w\rangle_{I_{-}}-\langle w\rangle_{I_{+}}}{\langle w\rangle_{I}}\right)^{2}
\langle v\rangle_{I}I\right)^{1/2}
$$
$$
 \leq C\|g\|_{2} \cdot \left(\sum|\langle fw^{1/2}\rangle_{I}|^{2} \left(
\frac{\langle w\rangle_{I_{-}}-\langle w\rangle_{I_{+}}}{\langle w\rangle_{I}}\right)^{2}
\langle v\rangle_{I}I\right)^{1/2}.
$$

Now we apply Lemma 2.1 to the second factor. Notice that the condition of Lemma
2.1 is satisfied (it is just (1.2)) we conclude that the second factor is at
most $C\|f\|_{2}$, and the theorem is proved for the case $v \in A_{\infty}$. The case $w \in
A_{\infty}$ is completely symmetric.

\bigskip

In the next theorem we require that 
\begin{equation}
w \in A_{\infty}(v).
\end{equation}

We remind the reader that there are many equivalent wordings of this assertion, which 
can be found in [St]. For us here it will be
important that this property means the following

$$\forall \e>0 \exists \delta>0, \, \forall I, E, E \subset I \,\,
\frac{v(E)}{v(I)} \leq \delta \Rightarrow \frac{w(E)}{w(I)} \leq \e.
$$

Remind also that this property is symmetric. 

\bigskip

{\bf Theorem 3.2} {\it If  $w$ is in $A_{\infty}(v)$ and if $w,v$ have the doubling property
then operator
$T_{0}$ is bounded if and only if (1.1) holds. In particular, in this situation
$\sup_{\e}\|T_{\e}\|_{L^2(w^{-1})\to L^2(v)} < \infty$ if and only if (1.1)--(1.3) hold.} 

\bigskip

Boundedness of $T_{0}$ is equivalent to (2.7). We have to be able to prove (2.7) starting
with our assumptions. An important particular case of (2.7) appears when one writes (2.7) for
$g = \chi_{J} v^{1/2}, \, f = \chi_{J} w^{1/2}$ where $J$ is an arbitrary dyadic interval. So
we have to be able to prove the following ``Carleson measure type inequality":

\begin{eqnarray*}
\forall J \in \mathcal{D}, \,\frac{1}{|J|}\sum_{I\subset J} \langle v\rangle_{I} \cdot \langle
w\rangle_{I}
\frac{|\langle v \rangle_{I_{-}}-\langle v\rangle_{I_{+}}|}{\langle v\rangle_{I}} \cdot
\frac{|\langle w \rangle_{I_{-}}-\langle w\rangle_{I_{+}}|}{\langle w\rangle_{I}} |I| 
\\[.2cm]
\leq C\sqrt{\langle v\rangle_J
\langle w\rangle_J}
\end{eqnarray*}.

It turns out that the proof of this inequality plays an important role in the proof of Theorem
3.2.

\bigskip

{\bf Lemma 3.3}  {\it If (1.1) holds then the above inequality  holds.
Moreover, for any
$\alpha
\in (0,1], \, \alpha
\neq\frac{1}{2},$ the following more general inequality holds}

\begin{eqnarray*}
\forall J \in \mathcal{D}, \,\frac{1}{|J|}\sum_{I\subset J} (\langle v\rangle_{I} \cdot \langle
w\rangle_{I})^{\alpha}
\frac{|\langle v \rangle_{I_{-}}-\langle v\rangle_{I_{+}}|}{\langle v\rangle_{I}} \cdot
\frac{|\langle w \rangle_{I_{-}}-\langle w\rangle_{I_{+}}|}{\langle w\rangle_{I}} |I| 
\\[.2cm]
\leq C(\alpha)({\langle v\rangle_J
\langle w\rangle_J})^{\min(\alpha, \frac{1}{2})}
\end{eqnarray*}.

We postpone the proof of Lemma 3.3 till the next section. We will see there that the constant
$C(\alpha)$ blows up when $\alpha$ approaches $\frac{1}{2}$. Now we use it to give the proof
of Theorem 3.2.

\bigskip

{\bf Proof of Theorem 3.2} Let $J$ be in $\mathcal{D}$ and let $\mathcal{I}$ be a disjoint family of
its dyadic subintervals having the following property

$$\frac{\langle F\rangle_{I}\langle G\rangle_{I}}{\langle w\rangle_{I}\langle v\rangle_{I}}
\geq B \frac{\langle F\rangle_{J}\langle G\rangle_{J}}{\langle w\rangle_{J}\langle
v\rangle_{J}}.
$$
Then $\frac{\cup_{I \in \mathcal{I}}v(I)}{v(J)} \leq \e$ and $\frac{\cup_{I \in
\mathcal{I}}w(I)}{w(J)} \leq \e$, where $\e$ depends on $B$ and is small if $B$ is large.

In fact, intervals $I$ are either of the type that
$\frac{\langle F\rangle_{I}}{\langle w\rangle_{I}} \geq \sqrt{B} \frac{\langle
F\rangle_{J}}{\langle w\rangle_{J}}$ or of the type that $\frac{\langle G\rangle_{I}}{\langle
v\rangle_{I}} \geq \sqrt{B} \frac{\langle G\rangle_{J}}{\langle v\rangle_{J}}$. Let
$\Omega_1$ be the union of the first type intervals, and let $\Omega_2$ denote the union of
the second type intervals, which are not in $\Omega_1$. Then obviously

$$w(\Omega_1)/w(J) \leq \delta, \, v(\Omega_1)/v(J) \leq \delta.
$$

We are ready to use that $w \in A_{\infty}(v)$ (which also means $v \in A_{\infty}(w)$)
to conclude that

\begin{equation}
\frac{\cup_{I \in \mathcal {I}}v(I)}{v(J)} \leq \e ,\, \frac{\cup_{I \in
\mathcal{I}}w(I)}{w(J)} \leq \e.
\end{equation}

To apply this remark let us fix $f,g$ from $L^2(\mathbb{T})$ and denote $F = fw^{1/2}, G =
gv^{1/2}$. Let $\langle F\rangle, \langle G\rangle$ are averages over $\mathbb{T}$.  Fix a
very large number
$C$ to be chosen later and let
${\mathcal J}_k$ denote maximal dyadic intervals such that

$$
\frac{\langle F\rangle_{I}\langle G\rangle_{I}}{\langle w\rangle_{I}\langle v\rangle_{I}}
\geq C^{k} \frac{\langle F\rangle\langle G\rangle}{\langle w\rangle\langle
v\rangle}.
$$

Also let ${\mathcal G}_k$ denote the collection of dyadic intervals, which are contained in some
$J$ from ${\mathcal J}_k$ and are not contained in any $J$ from ${\mathcal J}_{k+1}$.

Now we use Lemma 3.3 and (1.1) to conclude the following

\begin{eqnarray*}
\sum_{I\in {\mathcal G}_k } \langle F\rangle_{I} \langle
G\rangle_{I}
\frac{|\langle w \rangle_{I_{-}}-\langle w\rangle_{I_{+}}|}{\langle w\rangle_{I}} \cdot
\frac{|\langle v \rangle_{I_{-}}-\langle v\rangle_{I_{+}}|}{\langle v\rangle_{I}} |I|
\leq \\[.2cm]
\frac{\langle F\rangle\langle G\rangle}{\langle w\rangle\langle
v\rangle}C^{k+1} \sum_{I\in {\mathcal G}_k} \langle w\rangle_{I}\langle v\rangle_{I}\frac{|\langle w
\rangle_{I_{-}}-\langle w\rangle_{I_{+}}|}{\langle w\rangle_{I}} \cdot
\frac{|\langle v \rangle_{I_{-}}-\langle v\rangle_{I_{+}}|}{\langle v\rangle_{I}} |I|
\leq \\[.2cm]
\frac{\langle F\rangle\langle G\rangle}{\langle w\rangle\langle
v\rangle}C^{k+1}\sqrt{A}\sum_{J\in {\mathcal J}_k}\sqrt{\langle w\rangle_{J}\langle v \rangle_{J}} |J| 
\leq 
C\sqrt{A}\sum_{J\in {\mathcal J}_k}\frac{\langle F\rangle_{J}\langle G\rangle_{J}}{\langle
w\rangle_{J}\langle v\rangle_{J}}\sqrt{\langle w\rangle_{J}\langle v \rangle_{J}}|J|
=\\[.2cm]
C\sqrt{A}\sum_{J\in {\mathcal J}_k}\frac{\langle F\rangle_{J}\langle G\rangle_{J}}{\langle
w\rangle_{J}\langle v\rangle_{J}}\sqrt{w(J)v(J)}
\end{eqnarray*}

Here $A \df \sup_{J\in \mathcal{D}}\langle w\rangle_{J}\langle v\rangle_{J}$.

For an interval $J$ from ${\mathcal J}_k$ let us denote by $E_J$ the set $J\setminus \cup_{I \in
{\mathcal J}_{k+1}}I$. We would like to replace $w(J),v(J)$ in the last sum by $w(E_J),
v(E_J)$. Suppose for a moment that we can do that. Let us denote by $M_v$ the dyadic maximal
function with respect to measure $vdm$. Do similarly for $w$ to obtain $M_w$. Notice that
$\langle G\rangle_{J}/\langle v\rangle_J \leq M_v(gv^{-1/2})(x)$ for any $x \in J$.
Similarly, $\langle F\rangle_{J}/\langle w\rangle_J \leq M_w(fw^{-1/2})(x)$ for any $x \in
J$. Notice also that all $E_J$ are disjoint, and  we assume temporarily that, say, $w(J)
\leq 2w(E_J)$ and $v(J) \leq 2v(E_J)$.

Then we finish our estimate as follows:

\begin{eqnarray*}
\sum_{I\in \mathcal{D} } \langle F\rangle_{I} \langle
G\rangle_{I}
\frac{|\langle w \rangle_{I_{-}}-\langle w\rangle_{I_{+}}|}{\langle w\rangle_{I}} \cdot
\frac{|\langle v \rangle_{I_{-}}-\langle v\rangle_{I_{+}}|}{\langle v\rangle_{I}} |I|
\leq \\[.2cm]
\sum_{k}\sum_{I\in {\mathcal G}_k }...
\leq 
2C\sqrt{A}\sum_{k}\sum_{J\in {\mathcal J}_k}
\frac{\langle F\rangle_{J}\langle G\rangle_{J}}{\langle
w\rangle_{J}\langle v\rangle_{J}}\sqrt{w(E_J)v(E_J)}
\leq \\[.2cm]
2C\sqrt{A}(\int(M_w(fw^{-1/2})(x))^2 w(x)dx)^{1/2}(\int (M_v(gv^{-1/2})(x))^2
v(x)dx)^{1/2}
\leq \\[.2cm]
2C\sqrt{A}\|M_w(fw^{-1/2}\|_{L^2(w)}\|M_v(gv^{-1/2}\|_{L^2(v)}
\leq 
C^{'}\|f\|_{2} \|g\|_{2}
\end{eqnarray*}

We are left to choose $C$ so large that $w(J) \leq 2w(E_J)$ and $v(J) \leq 2v(E_J)$. To do
that let us notice that the maximality of $J$ in ${\mathcal J}_k$ and the doubling properties
imply that

$$
C^{k}\leq \frac{\langle F\rangle_{J}\langle G\rangle_{J}}{\langle w\rangle_{J}\langle
v\rangle_{J}} \leq DC^{k}.
$$

So if $I \subset J$ belongs to ${\mathcal J}_{k+1}$
we have that

$$ 
\frac{\langle F\rangle_{I}\langle G\rangle_{I}}{\langle w\rangle_{I}\langle
v\rangle_{I}}
\geq
\frac{C}{D}\frac{\langle F\rangle_{J}\langle G\rangle_{J}}{\langle w\rangle_{J}\langle
v\rangle_{J}}.
$$

The remark at the beginning of the proof of the theorem shows that if $C$ is much larger than
$D$ then the $v$-measure of $J\cap(\cup_{I\in {\mathcal J}_{k+1}}I)$ is smaller than any given
$\e$ (say, $\frac{1}{2}$). And the same is true about $w$-measure of this set. 

Thus $w(J) \leq 2w(E_J)$ and $v(J) \leq 2v(E_J)$, and the theorem is completely proved.

\setcounter{equation}{0}

\section{4. Bellman function and Carleson measures}

Here we are going to prove Lemma 3.3. We also apply the lemma to give another sufficient
condition for uniform boundedness of our $T_{\e}$.

In the proof of Lemma 3.3 the Bellman functions approach appears for the first time. It will
play the key role in the rest of the paper.

{\bf Proof of Lemma 3.3.} Without loss of
generality we assume that
\begin{equation}
\langle v\rangle_{I}\langle w\rangle_{I} \leq 1
\end{equation}

Let us first prove the case $\alpha \in (0,1/2)$.
Let us introduce the function
$$ 
\Phi(x,y) = \sup \frac{1}{|J|} \sum_{I\subseteq J} (\langle
v\rangle_{I}\langle w\rangle_{I})^{\alpha}
\left|
\frac{\langle v\rangle_{I_{-}}-\langle v\rangle_{I_{+}}}{\langle v\rangle_{I}}\right|
\left|\frac{\langle w\rangle_{I_{-}}-\langle w\rangle_{I_{+}}}{\langle w\rangle_{I}}\right|
|I|
$$ 
where supremum is taken over all nonnegative finite linear combinations $v, w$
of characteristic functions of dyadic intervals such that $v_{J} = x, w_{J} =
y$. Notice that $\Phi$ does not depend on $J$. Notice also that by the
definition $\Phi$ is concave. Moreover if $x =  \frac{x_{+}+x_{-}}{2}, \; y =
\frac{y_{+}+y_{-}}{2}$ then

\begin{equation}
\Phi(x,y) - \frac{(\Phi(x_{+},y_{+}) + \Phi(x_{-},y_{-}))}{2} \geq (xy)^{\alpha}
\left|\frac{x_{+}-x_{-}}{x}\right| \left|\frac{y_{+}-y_{-}}{y}\right|
\end{equation}

If only one could  prove that $\Phi(x,y) \leq C(\alpha)(xy)^{\alpha} \ldots \;$. Let us
reverse the argument. Suppose we can construct a concave $B$ in the domain
$ D = \{x \geq 0, y \geq 0, xy \leq 1\} $
such that

\begin{equation}
0 \leq B \leq (xy)^{\alpha}, \quad (x,y) \in \mathcal{D}
\end{equation}

and such that

\begin{eqnarray}
&& \hspace*{-.7in} B(x,y) - \frac{(B(x_{+},y_{+}) + B(x_{-},y_{-}))}{2} \geq
c_{\alpha} (xy)^{\alpha}
\left|\frac{x_{+}-x_{-}}{x}\right| \left|\frac{y_{+}-y_{-}}{y}\right| 
\end{eqnarray}
if $(x,y), (x_{+},y_{+}), (x_{-},y_{-}) \in D$ and if
$$ 
x = \frac{x_{+}+x_{-}}{2}, \quad y = \frac{y_{+}+y_{-}}{2}. 
$$
Then we are done.

And for $\alpha \in (0,1/2)$ it is easy to guess such a function:
$$ B(x,y) = (xy)^{\alpha}. $$

Let us consider it in a larger domain $ D_0 = \{x \geq 0, y \geq 0 \}$  and prove (4.4)
there.  Consider the function

$$ 
b(t) = B(x+t\xi, y+t\eta), \quad t \in [-1,1] 
$$
where $\xi = \frac{x_{+}-x_{-}}{2}, \; \eta = \frac{y_{+}-y_{-}}{2}$.
Clearly $b$
is concave, because $B$ is concave in $\{x \geq 0, y \geq 0\}$. We want to
estimate $b''(t)$ and to prove that

\begin{equation}
b''(t) \leq - c_{\alpha} \frac{|\xi\eta|}{x^{1-\alpha}y^{1-\alpha}}, \quad t \in [-1/2,1/2].
\end{equation}

If (4.5) were proved, then we could write
$$ b(0) - \frac{1}{2}(b(-1)+ b(1)) = -\int^{1}_{-1} (1-|t|)b''(t)dt \geq c_{\alpha}
\frac{|\xi\eta|}{x^{1-\alpha}y^{1-\alpha}} 
$$
which is exactly (4.4).

Let ${\mathcal H}_{x,y}(\xi,\eta)$ be the Hessian form of $B$ at point $(x,y)$
on the
vector $(\xi,\eta)$. Denote $x_{t} = x+t\xi, \; y_{t} = y+t\xi$. Then
$$ b''(t) = {\mathcal H}_{x_{t},y_{t}}(\xi,\eta). 
$$
Here is the computation of ${\mathcal H}_{x,y}(\xi,\eta)$ :

\begin{eqnarray*}
{\mathcal H}_{x_t,y_t}(\xi,\eta) = -\alpha x_t^{\alpha}
y_t^{\alpha}\left[(1-\alpha)(\frac{\xi}{x_t})^2 + (1-\alpha)(\frac{\eta}{y_t})^2 -
2\alpha\frac{\xi
\eta}{x_t y_t}\right] =\\[.2cm]
-\alpha x_t^{\alpha}
y_t^{\alpha}\left[(1-2\alpha)((\frac{\xi}{x_t})^2 + (\frac{\eta}{y_t})^2 ) + \alpha
(\frac{\xi}{x_t} +
\frac{\eta}{y_t})^2 \right] \leq \\[.2cm]
-(1-2\alpha)\alpha \frac{|\xi \eta|}{x_t^{1-\alpha}y_t^{1-\alpha}}
\end{eqnarray*}

Noticing that for $t \in [-1/2,1/2], \, x_t \geq \frac{x}{2}, \, y_t \geq \frac{y}{2}$ we get
(4.5),
which gives the proof of Lemma 3.3 for the case $\alpha < 1/2$.

Our key inequality (4.5) could be proved differently. Denote $t = xy, \; x_{\pm} = (1 \pm
\la)x,
\; y_{\pm} = (1\pm \mu)y,
\;
\la, \mu \in [-1,1]$ we come to proving that for $\la, \mu \in [-1,1]$ the
following holds
\begin{equation}
1 - \frac{1}{2}\left\{[(1-\la)(1-\mu)]^{\alpha} + [(1+\la)(1+\mu)]^{\alpha}\right\} \geq
c_{\alpha}|\la\mu|
\end{equation}
This elementary inequality is true when $\alpha \in (0,1/2)$, which can be checked by direct
computations.

Now let us consider the case $\alpha \in (1/2,1]$. Again we should
notice that it is enough to build a function $B$  such that
\begin{equation}
0 \leq B(x,y) \leq C(xy)^{1/2}, \quad (x,y) \in D,
\end{equation}
and
\begin{eqnarray}
&& \hspace*{-.7in} B(x,y) - \frac{B(x_{+},y_{+})+B(x_{-},y_{-})}{2} \geq
c_{\alpha}(xy)^{\alpha}
\left|
\frac{x_{+}-x_{-}}{x}\right| \left| \frac{y_{+}-y_{-}}{y}\right|
\end{eqnarray}
where both inequalities hold in $D = \{(x,y) : x \geq 0, y \geq 0, xy \leq
1\}$
and $x = \frac{x_{+}+x_{-}}{2}, \; y = \frac{y_{+}+y_{-}}{2}$.

Unfortunately, $B = (xy)^{1/2}$ does not work. It is not concave enough near the diagonal $x
= y$. Here is the function which satisfies (4.7), (4.8):
$$
B(x,y) = (xy)^{1/2} - \frac{1}{4} (xy)^{\alpha}. 
$$
Again there is nothing to prove in (4.7). It is a nonneggative function because  $xy \leq 1$
in $D$. 

Consider the function

$$ 
b(t) = B(x+t\xi, y+t\eta), \quad t \in [-1,1] 
$$
where $\xi = \frac{x_{+}-x_{-}}{2}, \; \eta = \frac{y_{+}-y_{-}}{2}$.
Clearly $b$
is concave, because $B$ is concave in $\{x \geq 0, y \geq 0\}$. We want to
estimate $b''(t)$ and to prove that
\begin{equation}
b''(t) \leq - c_{\alpha} \frac{|\xi\eta|}{x^{1-}y^{1-\alpha}}, \quad t \in [-1/2,1/2].
\end{equation}
If (4.9) were proved, then we could write
$$ b(0) - \frac{1}{2}(b(-1)+ b(1)) = -\int^{1}_{-1} (1-|t|)b''(t)dt \geq C(\alpha)
\frac{|\xi\eta|}{x^{1-\alpha}y^{1-\alpha}} $$
which is exactly (4.8).

Let ${\mathcal H}_{x,y}(\xi,\eta)$ be the Hessian form of $B$ at point $(x,y)$
on the
vector $(\xi,\eta)$. Denote $x_{t} = x+t\xi, \; y_{t} = y+t\xi$. Then
$$
b''(t) = {\mathcal H}_{x_{t},y_{t}}(\xi,\eta). 
$$

Here is the computation of ${\mathcal H}_{x,y}(\xi,\eta)$:

\pagebreak 

$$
{\mathcal H}_{x,y}(\xi,\eta) =
\frac{1}{4}(xy)^{1/2}[2(1-\alpha^{2}(xy)^{\alpha-1/2})
\frac{\xi\eta}{xy}-  
(1-\alpha(1-\alpha)(xy)^{\alpha-1/2})(\frac{\xi}{x})^{2}- 
$$
$$
(1-\alpha(1-\alpha)(xy)^{\alpha-1/2}) (\frac{\eta}{y})^{2}]
$$

Let us denote by $\rho$ the expression $\alpha^{2}(xy)^{\alpha-1/2}$ and by $r$ the
expression $\alpha(1-\alpha)(xy)^{\alpha-1/2}$. Then

\begin{eqnarray*}
\lefteqn{{\mathcal H}_{x,y}(\xi,\eta) = \frac{1}{4}(xy)^{1/2}\left\{-(1-\rho)
\left[\frac{\xi}{x} + \frac{\eta}{y}\right]^{2} \right.}
\\[.2cm]
&& \mbox{} \left. - (\rho-r) \left[\left(\frac{\xi}{x}\right)^{2} +
\left(\frac{\eta}{y}\right)^{2}\right]\right\} \leq -\frac{1}{4}
(\rho-r)(xy)^{1/2} \frac{|\xi\eta|}{xy}.
\end{eqnarray*}
Let us look at $\rho-r = (\alpha^{2}-\alpha(1-\alpha))(xy)^{\alpha-1/2} =
+\alpha(2\alpha-1)(xy)^{\alpha-1/2}$. It is positive, when $\alpha > 1/2$. In particular,

\begin{equation}
{\mathcal H}_{x_t,y_t}(\xi,\eta) \leq -c_{\alpha}
\frac{|\xi\eta|}{x^{1-\alpha}y^{1-\alpha}}
\end{equation}
if $t \in [-1/2,1/2]$ because $x_t \geq \frac{x}{2}, \, y_t \geq \frac{y}{2}$ for such $t$. 
Now (4.10) implies (4.9) immediately.  The proof of Lemma 3.3 is completed.

\bigskip

Let us show some corollaries of Lemma 3.3. Fix a
number $q \in (0,1)$. We denote by ${\mathcal F}_{k}, \; k = 0,1,2, \ldots$ the
family of all dyadic intervals $I$ such that
$$ 
q^{k+1} \leq v_{I}w_{I} \leq q^{k}.
$$
Consider the following measure in the disk
$$ 
\sigma_{k} = \sum_{I\in{\mathcal F}_{k}} \left|
\frac{v_{I_{-}}-v_{I_{+}}}{v_{I}}\right| \left|
\frac{w_{I_{-}}-w_{I_{+}}}{w_{I}}\right| |I| \delta_{c(I)}, 
$$
where $c(I)$ denotes the center of the box $Q(I)$ built over the dyadic arc $I$.

Let us remind that  measure $\sigma$ in $\mathbb{D}$ is called Carleson if $\sigma(Q(I))
\leq C|I|$ for all arcs $I$. The best constant in this inequality is called the Carleson
norm of
$\sigma$ and is denoted by $\|\sigma\|_{C}$.

\bigskip

{\bf Lemma 4.1.} {\it Measures $\sigma_{k}$ are Carleson measures and
$\|\sigma_{k}\|_{C} \leq B , \infty $.} 

\bigskip

The proof follows immediately from Lemma 3.3 with $\alpha \in (0,1/2)$. Let us
notice that in the case when $\langle v\rangle_{I}\langle w\rangle_{I}$ is bounded also from
below (so in the case when we can assume $w =
u^{-1} = v$, and we consider a one
weight problem) we get that $\sum^{\infty}_{k=0} \sigma_{k}$ is a Carleson measure.

\bigskip

We will see now how a slightly strengthen version of (1.1) becomes sufficient
for the  boundedness of $\|T_{0}\|_{L^2(w^{-1}) \to L^2(v)}$ (but it stops to be
necessary of course). We use the Carleson measure result that has been just proved.

\bigskip

{\bf Theorem 4.2.} {\it Suppose that} 

\begin{equation}
\forall \; I \in \mathcal{D}, \quad \langle v^{1+\eta}\rangle_{I} \langle
w^{1+\eta}\rangle_{I}
\leq A <
\infty
\end{equation}

{\it Then} $\|T_{0}\|_{L^{2}(w^{-1}) \to L^{2}(v)} \leq C(A) < \infty $. 

\bigskip

{\bf Remark.} There is a whole stream of results of this kind starting with
Fefferman-Phong theorem from [F] and continuing in [ChWW], [Z]. 

\bigskip

{\bf Proof.} Let us remind that the boundedness of $T_{0}$ is equivalent
to the following estimate for
all $L^{2}$-functions $f,g$

\begin{equation}
\sum_{I\in\mathcal{D}} \langle gv^{1/2}\rangle_{I} \langle fw^{1/2}\rangle_{I} \alpha_{I} \leq C
\|g\|_{2}
\|f\|_{2}.
\end{equation}
where

$$ 
\alpha_{I} = \left| \frac{\langle v\rangle_{I_{-}}-\langle
v\rangle_{I_{+}}}{\langle v\rangle_{I}}\right|
\left|
\frac{\langle w\rangle_{I_{-}}-\langle w\rangle_{I_{+}}}{\langle w\rangle_{I}}\right| |I|. 
$$
Clearly for $\gamma + \rho = 1/2, \gamma, \rho > 0$, we have

$$
\langle gv^{1/2}\rangle_{I} \leq \langle g^{p}\rangle_{I}^{1/p} \langle v^{\gamma
q_{1}}\rangle_{I}^{1/q_{1}}\langle v^{\rho q_{2}}\rangle_{I}^{1/q_{2}} 
$$
where $1/ p + 1/q_{1} + 1/q_{2} = 1$. Choose $\gamma = \rho = 1/4, \; q_{1} = 4,
\; q_{2} = 4(1+\eta)$ to have $p = 4 \frac{1+\eta}{2+3\eta} < 2$. 
Then

\begin{eqnarray*}
\langle gv^{1/2}\rangle_{I} & \leq & \langle g^{p}\rangle_{I}^{1/p}
\langle v^{1+\eta}\rangle_{I}^{1/q_{2}}\langle v\rangle_{I}^{1/4}, \\[.2cm]
\langle fw^{1/2}\rangle_{I} & \leq &
\langle f^{p}\rangle_{I}^{1/p}\langle w^{1+\eta}\rangle_{I}^{1/q_{2}}\langle
w\rangle_{I}^{1/4}.
\end{eqnarray*}

Remind that measures from Lemma 4.1
$$ 
\sigma_{k} = \sum_{I\in{\mathcal F}_{k}} \alpha_{I} \delta_{c(I)} 
$$
are proved to be uniformly Carleson. Thus the measure
$$ 
\sigma \stackrel{def}{=} \sum_{I\in\mathcal{D}} (v_{I}w_{I})^{1/4} \alpha_{I}\delta_{c(I)}
\leq \sum_{k} (q^{k})^{1/4} \sigma_{k} 
$$
is a Carleson measure. Let us notice that $2/p > 1$, and let us use the Carleson embedding
theorem in
$L^{2/p}$ (see [G]) to estimate the sum in (4.12).

\begin{eqnarray*}
\lefteqn{ \sum_{I\in\mathcal{D}} \langle gv^{1/2}\rangle_{I}\langle fw^{1/2}\rangle_{I}\alpha_{I}
\leq
\sup_{I}(\langle v^{1+\eta}\rangle_{I}\langle w^{1+\eta}\rangle_{I})^{1/q_{2}}} \\[.2cm]
&& \cdot \sum_{I\in\mathcal{D}}
\langle g^{p}\rangle_{I}^{1/p}\langle f^{p}\rangle_{I}^{1/p}(\langle
v\rangle_{I}\langle w\rangle_{I})^{1/4}\alpha_{I}
\\[.2cm] && \leq \left(\sum_{I\in\mathcal{D}}
\langle f^{p}\rangle_{I}^{2/p}(
\langle v\rangle_{I}\langle w\rangle_{I})^{1/4}\alpha_{I}\right)^{1/2} \left(
\sum_{I\in\mathcal{D}}\langle
f^{p}\rangle_{I}^{2/p}(\langle v\rangle_{I}\langle
w\rangle_{I})^{1/4}\alpha_{I}\right)^{1/2}
\\[.2cm] && \leq C(p)\|\sigma\|_{C} \|g^{p}\|^{1/p}_{L^{2/p}} \|f^{p}\|^{1/p}_{L^{2/p}}
\leq C \|f\|_{L^{2}} \|g\|_{L^{2}}.
\end{eqnarray*}

We are done. 

\setcounter{equation}{0}

\section{5. The necessary and sufficient conditions: $\!$Bilinear weighted imbedding
theorem}

We saw in Section 2 that given necessary conditions (1.1), (1.2), (1.3) the uniform
boundedness of $\|T_{\e}\|_{L^2(w^{-1})} \to L^2(v)$ is equivalent to the following estimate
for all $L^{2}$-functions $f,g$

\begin{equation}
\sum \langle gv^{1/2}\rangle_{I} \langle fw^{1/2}\rangle_{I} \alpha_{I} \leq C\|g\|_{2}
\|f\|_{2}
\end{equation}
where  $\alpha_{I} = \left|\frac{\langle v\rangle_{I_{-}}-\langle
v\rangle_{I_{+}}}{\langle v\rangle_{I}}\right|
\left|\frac{\langle w\rangle_{I_{-}}-\langle w\rangle_{I_{+}}}{\langle w\rangle_{I}}\right|
|I|$. We may (and should) consider (5.1) as a bilinear weighted imbedding estimate.

Consider the kernel function
$$
k(x,y) = \sum_{I\in\mathcal{D}} \frac{1}{|I|^{2}} \chi_{I}(x) \chi_{I}(y) \alpha_{I}
$$
It is obvious that (5.1) means the estimate for any nonnegative measurable $G$:

\begin{equation}
\int \left(\int k(x,y) G(y)w(y)dy\right)^{2} v(x)dx \leq C \int G^{2} w \, dy
\end{equation}

The inequalities (5.2) have been extensively studied in the works of  Kalton, Sawyer,
Verbitsky, and Wheeden.  The next result is in the vein of
those works. Inequalities (5.2) can be most probably proved using the combination of ideas
frim the papers of Kalton, Verbitsky [KV] and Sawyer, Wheeden [SW]. Following [KV] we can
introduce the new metric $d(x,y) := \frac{1}{k(x,y)}$, where $k$ is the kernel of $T_0$
and was written above. It is clearly a metric, and balls in this metric are just all
dyadic intervals. Then kernel $k$ satisfies the regularity condition from [SW]. Using
the main result of [SW] we could have given an alternative proof of  inequalities (5.2).
Unfortunately it is not clear why the new metric space has certain regularity properties.
For example, in [KV] one requires the property that all annuli be nonempty. This is false
in our new metric space. However, it is not clear to us how essential are those
regularity properties of the metric space for the application of the technique of [SW].

So we will give another proof.

Plug into (5.2) the charecteristic functions of the intervals. Then we get a condition
necessary  for (5.2). Do the same with the dual inequality. Then we get another necessary
condition. These are are the conditions of the Sawyer type:

\begin{eqnarray}
\forall \; J \in \mathcal{D}, && \frac{1}{|J|}\int_{J} \left(\int_{J} k(x,y)w(y)dy\right)^{2}
v(x)dx \leq C\langle w\rangle_J, \\[.2cm]
\forall \; J \in \mathcal{D}, && \frac{1}{|J|}\int_{J} \left( \int_{J} k(x,y)v(x)dx\right)^{2}
w(y)dy \leq C\langle v\rangle_J.
\end{eqnarray}

In other words, for all dyadic arcs $J$

\begin{equation}
\frac{1}{|J|}\int_{J} \left(\sum_{I\subset J} \frac{1}{I} \chi_{I} \langle w\rangle_{I}
\alpha_{I}\right)^{2} v
\, dx \leq C\langle w\rangle_J,
\end{equation}
\begin{equation}
\frac{1}{|J|}\int_{J} \left( \sum_{I\subset J} \frac{1}{I}
\chi_{I}\langle v\rangle_{I}\alpha_{I}\right)^{2} w
\, dx \leq C\langle v\rangle_J.
\end{equation}

Let us denote by $T_{0}$ the operator with kernel $k$.

\bigskip

{\bf Theorem 5.1.} {\it Let $\{\alpha_{I}\}$ be a nonnegative sequence. The
following assertions are equivalent}

\begin{enumerate}
\item[1)] {\it Operator $T_{0}$ is bounded from $L^{2}(w^{-1})$ to
$L^{2}(v)$;}
\item[2)] {\it Estimate (5.1) holds;}
\item[3)] {\it Estimate (5.2) holds;}
\item[4)] {\it For any $J \in \mathcal{D}$ (5.5) and (5.6) hold with the same
constant $C$ independent of $J$.} 
\end{enumerate}

\bigskip

{\bf Remark:} Clearly (5.5) (or (5.6) as well) implies

\begin{equation}
\frac{1}{|J|}\sum_{I\subset J}\langle v\rangle_{I}\langle w\rangle_{I}\alpha_{I} \leq C
\sqrt{\langle v\rangle_J \langle w\rangle_J},
\end{equation}
which is the conclusion of Lemma 3.3 for the case $\alpha = 1$, and so (5.7)
follows from just (1.1). Unfortunately, (5.5), (5.6) themselves do not follow
from (1.1).

On the other hand (5.7) is a particular case of our key estimate (5.1). In fact,
(5.1) becomes (5.7) if $g = \chi_{J}v^{1/2}, f = \chi_{J}w^{1/2}$. The
conclusion is that the particular case of (5.1) for test functions
$\chi_{J}v^{1/2},\chi_{J}w^{1/2}$ follows from (1.1). But this is not the case for
the full of (5.1). However, (5.5) and (5.6) are in fact the conditions obtained by testing
(5.1) by certain test functions. We are going to repeat what has been already said above.
Namely, it
is easy to see that (5.5) is the testing of (5.1) by
$f = \chi_{J}w^{1/2}$ and $g = g^{J}$, where $g^{J}$ is any $L^{2}$-function
supported on $J$. Similarly (5.6) is the testing of (5.1) by $g =
\chi_{J}v^{1/2}$ and $f = f ^{J}$, where $f^{J}$ is any $L^{2}$-function
suppoprted on $J$.

Let us conclude that in a particular case $v = w$ the bilinear weighted
imbedding estimate (5.1) becomes a usual weighted $L^{2}$-imbedding estimate for
linear operator $f \to \{(fw^{1/2})_{I}\}_{I \in \mathcal{D}}$:

\begin{equation}
\sum \langle fw^{1/2}\rangle_{I}^{2} \alpha_{I} \leq C\|f\|^{2}_{2}.
\end{equation}

Test conditions (5.5), (5.6) become

\begin{equation}
\forall \; J \in \mathcal{D}, \; \frac{1}{|J|} \sum_{I \subset J} \langle w\rangle_{I}^{2}
\alpha_{I}
\leq C\langle w\rangle_{J}.
\end{equation}

Let us start our proof of Theorem 5.1 with a special proof of the classical
equivalence (5.8) $\Leftrightarrow$ (5.9) (see [S1], [S2], [TV]).

\setcounter{equation}{0}

\section{6. The necessary and sufficient conditions: the proof of Theorem
5.1 for the case $w =v$}

{\bf Theorem 6.1} (5.8) $\Leftrightarrow$ (5.9). 

\bigskip

{\bf Proof.} One needs to prove only the implication (5.9) $\Rightarrow$ (5.8).
Let us consider the function
$$ 
\Phi_{J}(X,x,w) = \sup \frac{1}{|J|} \sum_{I\subseteq J} \langle fw^{1/2}\rangle_{I}^{2}
\alpha_{I} 
$$
where supremum is taken over $f \geq 0, \; \langle fw^{1/2}\rangle_{J} = x, \;
\langle f\rangle_{J}^{2} = X,
\; \langle w\rangle_{J} = w$. Clearly, this restricts us to
$$ 
D = \{(X,x,w) : x^{2} \leq X w, \; X, x, w \geq 0\}. 
$$
The family of functions $\Phi_{J}$ is concave in the following sense.
Let  $X = \frac{X_{-} + X_{+}}{2}, \, x = \frac{x_{-} + x_{+}}{2}, \, w = \frac{w_{-} +
w_{+}}{2}$. Then

\begin{equation}
\Phi_{J}(X,x,w) - \frac{1}{2}(\Phi_{J_{-}}(X_{-},x_{-},w_{-}) +
\Phi_{J_{+}}(X_{+},x_{+},w_{+}) \geq \frac{1}{|J|} x^{2}\alpha_{J}
\end{equation}

In fact, fixing averages on $J_{-}, J_{+}$ separately leads to a smaller set of functions
when if fixing the averages only on $J$ (we keep the compatability by saying that the
averages over $J$ are arithmetic means of averages over $J_{-}, J_{+}$. 

One wishes to prove that
$\Phi_{J} \leq CX$. 
This certainly would give the result. But this is certainly false because
we did not use the condition (5.9) in the definition of $\Phi_{J}$.

Let us try to correct this by introducing
$$ 
\Psi_{J}(X,x,w) = \sup \frac{1}{|J|} \sum_{I \subseteq J} \langle fw^{1/2}\rangle_{I}^{2}
\alpha_{I}
$$
where supremum is taken over the same set of conditions including $\langle w\rangle_{J} = w$
plus the condition that $\{\langle w\rangle_{I}\}_{I \subset J}$ are so distributed that

\begin{equation}
\frac{1}{|J'|} \sum_{I\subseteq J'} \langle w\rangle_{I}^{2} \alpha_{I} \leq C\langle
w\rangle_{J'},
\quad
\forall
\; J'
\subset J
\end{equation}

Unfortunately, it is now not clear why (6.1) would hold with $\Psi_{J}$
replacing $\Phi_{J}$. Conditioning as (6.1) makes unclear why $\{\Phi_{J}\}$ has
even simple concavity property. In fact, by passing to $J_{-}, J_{+}$ from $J$ we did not
deminished the set over which the supremum is taking. Unlike the case of $\Phi$ above we
actually relaxed our requirements. Notice that inequality (6.2) for $J^{'}=J$ will not be
required anymore after passing to  $J_{-}, J_{+}$.

But one can overcome this difficulty by fixing the ``main" sum in (6.1) and making its value
$M$ the new variable. 
To do that let us consider 
$$ 
\Phi_{J}(X,x,w,M) \stackrel{def}{=} \sup \frac{1}{|J|} \sum_{I\subseteq J} \langle
fw^{1/2}\rangle_{I}^{2}
\alpha_{I} 
$$
where supremum is taken over such functions that
$f \geq 0, \; \langle fw^{1/2}\rangle_{J} =
x,
\;
\langle f\rangle_{J}^{2} = X,
\; \langle w\rangle_{J} = w$, and over sequences $\{\alpha_{I}\}$ and functions such that

\begin{equation}
\frac{1}{|J'|} \sum_{I\subseteq J'} \langle w\rangle_{I}^{2} \alpha_{I} \leq C\langle
w\rangle_{J'},
\quad
\forall
\; J'
\subset J,
\end{equation} 
and such that

\begin{equation}
 M = \frac{1}{|J|} \sum_{I\subseteq J} \langle w\rangle_{I}^{2} \alpha_{I}. 
\end{equation}

Notice that thus defined $B$ does not depend on $J$, and that it is cleary defined and concave
in the domain

\begin{eqnarray*}
D_{B} & = & \{(X,x,w,M) : x^{2} \leq Xw, \; (X,x,w,M) \geq 0, \; M
\leq Cw \}.
\end{eqnarray*}

On the top of concavity we know the following. Let $h = M -\frac{M_{-} +M_{+}}{2}$. Clearly,

\begin{equation}
B(X,x,w,M) -\frac{B(X,x,w,M_{-})+B(X,x,w,M_{+})}{2} \geq \frac{x^2}{w^2}h.
\end{equation}

Thus we conclude that $B$ satisfies the following properties

\begin{eqnarray}
0 & \leq & B \leq X, \quad \mbox{on} \quad \mathcal{D}_{B}, \\[.2cm]
\frac{\partial B}{\partial M} & \geq & \frac{x^{2}}{w^{2}}, \quad \mbox{on}
\quad \mathcal{D}_{B}, \\[.2cm]
d^{2}B & \leq & 0, \quad \mbox{on} \quad \mathcal{D}_{B},
\end{eqnarray}

Instead of looking for the exact formula for $B$ let us consider
$$ 
B(X,x,w,M) = CX- C\frac{x^{2}}{w+M}, 
$$
which satisfies all the properties (6.8)-(6.10). To finish the proof of (5.8) it is enough to
consider triples $ x=\langle fw^{1/2}\rangle_{J},
\;
X=\langle f\rangle_{J}^{2} ,
\;w = \langle w\rangle_{J}, x_{\pm}=\langle fw^{1/2}\rangle_{J_{\pm}},
\;
X_{\pm}=\langle f\rangle_{J_{\pm}}^{2} ,
\;w_{\pm} = \langle w\rangle_{J_{\pm}}$ and $M,M_{\pm}$ such that $M -\frac{M_{-} +M_{+}}{2}
= \frac{1}{|J|}w^2 \alpha_{J}$, and to prove that
$$
B(X,x,w,M) -\frac{B(X_{-},x_{-},w_{-},M_{-})+B(X_{+},x_{+},w_{+},M_{+})}{2} \geq \frac{1}{|J|}
x^2
\alpha_{J}.
$$
But the last inequality follows immediately from (6.7) and the concavity of $B$.

{\bf Remark.} It is interesting to notice that one can give another euristic explanation of
appearance of new variable $M$. Let us think that all averages are variables and all sums are
functions defined on the set of all averages involved. Then $M$ is a pretty concave function
of
$w$ in the following sense:

\begin{equation}
M(w_{J}) - \frac{1}{2}[M(w_{J_{+}}) + M(w_{J_{-}})] \geq \frac{1}{|J|} w^{2}_{J}
\alpha_{J}.
\end{equation}
Thus it is resonable to look for $Q(X,x,w)$ in the form $Q = B(X,x,w,M(w))$, where $B$ has
a large derivative with respect to $M$. Then we can use the following formula for the
Hessian of the composition

\begin{equation}
d^2 Q = d^2 B + \frac{\partial{B}}{\partial{M}}\cdot d^2 M.
\end{equation}

\setcounter{equation}{0}

\section{7. Necessary and sufficient conditions. The proof of the bilinear weighted
imbedding theorem}

We try to repeat the argument of the previous section. First of all it would be nice to
discern the variables of the future Bellman function. Some of them are ready right away. These
are $X=\langle f^2\rangle_J, x=\langle fw^{1/2}\rangle_J, w =w_J, Y=\langle g^2\rangle_J,
y=\langle gv^{1/2}\rangle_J, v =v_J$. Imitating the considerations of the previous section
we also fix $M  =  \frac{1}{|J|} \int_{J} \left(\sum_{I\subseteq J}
\frac{1}{|I|}
\chi_{I} w_{I}\alpha_{I}\right)^{2}v \, dx$ and $N  =  \frac{1}{|J|} \int_{J} \left(
\sum_{I\subseteq J}
\frac{1}{|I|} \chi_{I} v_{I}\alpha_{I}\right)^{2} w \, dx$. By doing this we are making  the
following function concave:
$$
\Phi(X,x,w,Y,y,v,M,N) \df \sup\frac{1}{|J|}\sum_{I\subset J}\langle
fw^{1/2}\rangle_{I}\langle gv^{1/2}\rangle_{I} \alpha_{I}
$$
where the supremum is taken over nonnegative functions with ``main" averages fixed as above
and such that

\begin{eqnarray*}
\forall I\in \mathcal{D}, \, I\subset J  \,
\frac{1}{|I|} \int_{I} \left(\sum_{I'\subseteq I}
\frac{1}{|I'|}
\chi_{I'} x_{I'}\alpha_{I'}\right)^{2}v \, dx \leq C w_{I}  \\[.2cm]
\forall I\in \mathcal{D}, \, I\subset J  \,
\frac{1}{|I|} \int_{I} \left( \sum_{I'\subseteq I}
\frac{1}{|I'|} \chi_{I'} v_{I'}\alpha_{I'}\right)^{2} w \, dx \leq C v_{I}
\end{eqnarray*}

We could try to use $\Phi$ as Bellman function of our problem having hope that it gives our
main inequality (5.1). Actually it does not. It gives a certain inequality, but it is weaker
than the one we need. But on the other hand we used all possible variables (averages,
conditions). There are no more averages involved, and there are no more conditions (if we
believe that our bilinear imbedding theorem is true).

First let us find the way out of this impass by euristic considerations. Notice that in the
previous sections once the Bellman function is found we apply a certain inequality not to
the Bellman function itself but rather to its composition with martingales. In fact, all
these
$\langle f^2\rangle_J,...,v =\langle v\rangle_J$ are martingales. There are also
supermartingales. For example,  
$$
M_{J}  =  \frac{1}{|J|} \int_{J} \left(\sum_{I\subseteq J}
\frac{1}{|I|}
\chi_{I} w_{I}\alpha_{I}\right)^{2}v \, dx
$$
is a supermartingale. Let us pay attention to its
``discrete Laplacian" $M_{J} - \frac{M_{J_{-}} + M_{J_{+}}}{2}$. If we consider a
supermartingale $m_{J} \df 
\frac{1}{|J|} \sum_{I\subseteq J} \langle w\rangle_{I}^{2} \alpha_{I}$ from the
previous section we see that its ``discrete Laplacian" $m_{J} - \frac{m_{J_{-}} +
m_{J_{+}}}{2}$ is equal to $\frac{1}{|J|}\langle w\rangle_{J}^2 \alpha_{J}$. We notice that
it involves only martingale $\langle w\rangle_{J}$, which has been already chosen as a
variable (called
$w$) of the Bellman function.

What happens with ``discrete Laplacian" of supermartingale $M_{J}$? When we calculate it, we
get

\begin{equation}
M_{J} - \frac{M_{J_{-}} + M_{J_{+}}}{2} = \frac{1}{|J|}\langle w\rangle_{J}K_{J}\alpha_{J},
\end{equation}
where 
$$
\K_{J} \df \frac{1}{|J|}\sum_{I\subset J}\langle w\rangle_{I}\langle v\rangle_{I} \alpha_{I}.
$$
So the ``second derivative" of supermartingale $M_{J}$ involves not only the martingale
$\langle w\rangle_{J}$ but also a new supermartingale $K_{J}$.

The supermartingale $K_{J}$, which appeared as ``discrete Laplacian" of our ``natural
 variable" supermartingale $M$ also should be embodied by a variable--naturally let us call it
$K$.

One may now try to continue this process by taking the ``second derivative" of 
$K_{J}$:
 
\begin{equation}
K_{J} - \frac{K_{J_{-}} + K_{J_{+}}}{2} = \frac{1}{|J|}\langle w\rangle_{J}\langle
v\rangle_{J} \alpha_{J}
\end{equation}
to make sure that it involves only already considered martingales
$\langle w\rangle_{J}, \langle v\rangle_{J}$. Fortunately no new variables appear after
$K$.
 
{\bf The recipe}. Let us repeat that we are calculating the ``discrete Laplacian" of the
composition of a certain (concave) function (Bellman function) and several martingales and
supermartingales. It is only too natural that in this calculation the ``discrete Laplacians"
of all these martingales and supermartingales appear. For martingales they will be zero, for
supermartingales they are positive and involve our martingales and may be some new
(super)martingales. If this happens, this new (super)martingales should be added to the list
of (super)martingales, which must become the variables. 

In the situation we have now we have another supermartingale $N_{J}$ to which we have to
apply this algorithm. It is given by the formula:
$$
N_{J}  =  \frac{1}{|J|} \int_{J} \left(\sum_{I\subseteq J}
\frac{1}{|I|}
\chi_{I} v_{I}\alpha_{I}\right)^{2}w \, dx.
$$

And so its ``discrete laplacian" can be calculated as follows:

\begin{equation}
N_{J} - \frac{N_{J_{-}} + N_{J_{+}}}{2} = \frac{1}{|J|}\langle v\rangle_{J}K_{J}\alpha_{J},
\end{equation}
where $K_{J}$ is the same supermartingale we found by calculating the ``discrete Laplacian"
of $M_{J}$ above.

Now we have all the variables: $X,x,w,Y,y,v,K,M,N$. Let us make the notation for this
$9$-tuple and also for $X,x,w,Y,y,v,M,N$ and $X,x,w,Y,y,v,K$.
$$
a \df (X,x,w,Y,y,v,K,M,N), \, b \df (X,x,w,Y,y,v,M,N),\, c \df (X,x,w,Y,y,v,K).
$$
 
For technical reason we will be looking for $B$ in the form
$B(a) = Q(b) + P(c)$.

The domain of definition is 
$$
D = D_{B}  =\{ a = (X,x,w,Y,y,v,K,M,N): x^2 \leq Xw, y^2 \leq Yv,
$$
$$
M \leq Cw, N \leq Cv, K \leq C\sqrt{wv}, a \geq 0 \}.
$$

The only thing requiring the clarification is why $ K \leq C\sqrt{wv}$  here?
But $K$ stands for supermartingale $K_{J}$ for which this ineqality was proved in Lemma 3.3.
So this ineqality represents a natural restriction on $K$. 

What inequalities on $B$ would be sufficient for us? It should be the inequality on the
composition of
$B$ with all our (super)martingales. So consider $a_{J} =
(\langle f^2\rangle_{J},\langle fw^{1/2}\rangle_{J},w_{J},\langle g^2\rangle_{J},\langle
gv^{1/2}\rangle_{J},v_{J},K_{J},M_{J},N_{J})$ and $a_{J_{\pm}}$ denoting the same thing for
$J_{\pm}$ instead of $J$.
Here are the inequality we want to obtain

\begin{equation}
B(a_{J}) - \frac{B(a_{J_{-}})+B(a_{J_{+}})}{2} \geq
c\frac{1}{|J|}\langle fw^{1/2}\rangle_{J}\langle gv^{1/2}\rangle_{J}\alpha_{J}
\end{equation}

\begin{equation}
0 \leq B(a) \leq C(X+Y)
\end{equation}

If we have both of these inequalities then moving from $J$ to $J_{\pm}$ and then continuing
these process to ``sons" of $J_{\pm}$ et cetera... we obtain the inequality

\begin{equation}
\frac{1}{|J|}\sum_{I\subset J}\langle
fw^{1/2}\rangle_{I}\langle gv^{1/2}\rangle_{I} \alpha_{I} \leq C(\langle
f^2\rangle_{J}+\langle g^2\rangle_{J})
\end{equation}

The use of homogenuity of the left part (it does not change under the replacement of $f$ by
$tf$, $g$ by $t^{-1}g$) would finish the proof of our bilinear weighted imbedding theorem
because we obtain

\begin{equation*}
\frac{1}{|J|}\sum_{I\subset J}\langle
fw^{1/2}\rangle_{I}\langle gv^{1/2}\rangle_{I} \alpha_{I} \leq C\sqrt{\langle
f^2\rangle_{J}\langle g^2\rangle_{J}} = C\|f\|_{L^2_J}\|g\|_{L^2_J}.
\end{equation*}

\bigskip 

{\bf Lemma 7.1} {\it To have (7.4) it is sufficient to 
have the following inequalities}
\begin{equation}
B(a) -\frac{B(a_{-})+B(a_{+})}{2} \geq \gamma_1\frac{xy}{wv}(K - \frac{K_{-}+K_{+}}{2})
\end{equation}
{\it if} 
\begin{equation}
\frac{x^2}{w}K + \frac{y^2}{v}K \leq Cxy
\end{equation}
{\it and}
\begin{equation}
B(a) -\frac{B(a_{-})+B(a_{+})}{2} \geq
\gamma_2 \frac
{x^2}{w^2}(M - \frac{M_{-}+M_{+}}{2})+\gamma_3 \frac{y^2}{v^2}(N -
\frac{N_{-}+N_{+}}{2}),
\end{equation}
{\it if}
\begin{equation}
\frac{x^2}{w}K + \frac{y^2}{v}K \geq Cxy
\end{equation}
{\it Here} $a \geq \frac{a_{-}+a_{+}}{2}$. {\it  And } $\gamma_{i}, \, i =1,2,3$, {\it are
positive constants.}

\bigskip

{\bf Proof}. Let us calculate the right parts of (7.7) and (7.9) for
 
$$
a = (X,x,w,Y,y,v,K,M,N)=
(\langle f^2\rangle_{J},\langle fw^{1/2}\rangle_{J},w_{J},
$$
$$
\langle g^2\rangle_{J},\langle
gv^{1/2}\rangle_{J},v_{J},K_{J},M_{J},N_{J}) ,
$$
and $a_{\pm}$ 
having the same meaning with $J$
replaced by $J_{\pm}$. We do this calculation by using (7.1)-(7.3). In both cases, the right
part is bigger than

\begin{eqnarray*}
\frac{1}{|J|} xy \alpha_{J} = \frac{1}{|J|}\langle fw^{1/2}\rangle_{J}\langle
g^2\rangle_{J}\alpha_{J}.
\end{eqnarray*}

Notice that we use (7.10) when we calculate the right part of (7.9). At any rate in both
cases we get (7.4).

\bigskip

Now we are left to find $B$ satisfying (7.7) and (7.9).
As we mentioned we will be looking for $B$ in the form $B(a)=Q(b)+P(c)$ where $
a = (X,x,w,Y,y,v,K,M,N), \, b = (X,x,w,Y,y,v,M,N),\, c=(X,Y,x,w,y,v,K)$.

\bigskip

{\bf Lemma 7.2} {\it Function} $Q = CX+CY - \frac{x^2}{w + M} - \frac{y^2}{v + N}$ {\it
satisfies the inequality}

\bigskip

\begin{equation}
Q(b) - \frac{Q(b_{-})+Q(b_{+})}{2} \geq
\gamma_2 \frac
{x^2}{w^2}(M - \frac{M_{-}+M_{+}}{2})+\gamma_3 \frac{y^2}{v^2}(N -
\frac{N_{-}+N_{+}}{2}),
\end{equation}
{\it if} $ b \in D_{Q} \df \{ b = (X,x,w,Y,y,v,M,N): x^2 \leq Xw, y^2 \leq Yv,\\[.2cm]
M \leq Cw, N \leq Cv, b \geq 0 \}$.

\bigskip

{\it Proof}. This is a direct computation. It easily follows from the following infinitesimal
estimates for $Q$:

\begin{eqnarray*}
\frac{\partial Q}{\partial M} & \geq & \frac{cx^{2}}{w^{2}}, \,\\[.2cm] 
\frac{\partial Q}{\partial N} & \geq & \frac{cy^{2}}{v^{2}}, \,\\[.2cm] 
d^{2}Q & \leq & 0.
\end{eqnarray*}

\bigskip

{\bf Lemma 7.3} {\it Function} 
$$
P = CX+CY -  \sup_{0 < s < \infty} \left(\frac{x^{2}}{w+sK} +
\frac{y^{2}}{v+s^{-1}K}\right)
$$
{\it
satisfies the inequality}

\begin{equation}
P(c) - \frac{P(c_{-})+P(c_{+})}{2}  \geq \gamma_1\frac{xy}{wv}(K -
\frac{K_{-}+K_{+}}{2})
\end{equation}
{\it if} 
\begin{equation*}
\frac{x^2}{w}K + \frac{y^2}{v}K \leq Cxy
\end{equation*}
{\it if} $ c \in D_{P} \df \{ b = (X,x,w,Y,y,v,K): x^2 \leq Xw, y^2 \leq Yv,\\[.2cm]
K \leq C\sqrt{wv}, c \geq 0 \}$.

\bigskip

{\it Proof}. We hope that the formula for this function is as much surprising for the reader
as it was for us. But once the formula is written the rest is a direct computation. It easily
follows from the following infinitesimal estimates for $P$.

Optimal $s$ can be found from the equation

\begin{equation}
\frac{sK}{(w+sK)^{2}} x^{2} = \frac{s^{-1}K}{(v+s^{-1}K)^{2}} y^{2}
\end{equation}

Now we can make the direct computation of $\frac{\partial P}{\partial K}$. If $sK \leq cw$
and
$s^{-1}K \leq Cv$, then
$$ 
\frac{\partial P}{\partial K} \asymp \frac{x^{2}}{w^{2}} s +
\frac{y^{2}}{v^{2}} s^{-1} \geq \frac{xy}{wv}. 
$$
It is left to show that if $sK >> w$ or if $s^{-1}K >> v$ (the symbol $>>$
means ``much larger than'') then we are automatically in the area  where
$$ 
\left( \frac{x^{2}}{w} + \frac{y^{2}}{v}\right)K \geq c \, xy. 
$$
In fact, if $s^{-1}K >> v$ then $sK << w$. This is because we are in the domain where
$K \leq C\sqrt{wv}$
Then (7.12) becomes
$$ 
\frac{sK}{w^{2}} x^{2} = \frac{1}{s^{-1}K} y^{2} 
$$

Thus $\frac{x^{2}}{y^{2}} \asymp \frac{w^{2}}{K^{2}}$ and $\frac{x^{2}}{w}
K \geq cxy \frac{w}{K} \frac{K}{w} = cxy$. Similarly, if $sK >> w$
and thus $s^{-1}K << v$, equality (7.12) gives $\frac{y^{2}}{x^{2}} \asymp
\frac{v^{2}}{K^{2}}$ and so $\frac{y^{2}}{v} K \geq cxy
\frac{v}{K} \frac{K}{v} = cxy$. Lemma 7.3 is completely proved.

So bilinear weighted imbedding theorem is fully proved.

\bigskip

\setcounter{equation}{0}

\section{8. Hilbert transform. Sufficient conditions via Green's potentials.}

Let $H$ stands for the Hilbert transform on the circle $\mathbb{T}$. In other words, operator
$H$ (defined first on trigonometric polynomials) acts by the formula:

$$
H(\sum a_k e^{i\theta k}) \stackrel{def}{=} -i\sum_{k \geq 0} a_k e^{i\theta k} + i\sum_{k <
0} a_k e^{i\theta k} \,.
$$

Throughout this section we use the notation $f(z)$ for the Poisson extension of the function
$f \in L^{1}(\mathbb{T})$ to the disc $\mathbb{D}$ evaluated at the point $ z \in \mathbb{D}$.
So, for example,
$f(z)^2$ and $f^2(z)$ are different in general.
 
The two weights
estimates for the Hilbert transform appear naturally in the theory of Hankel and Toeplitz
operators and in the perturbation theory of linear operators (including the perturbation of
differential operators). 

For the case of the Hilbert transform there exists the approach  of
 M.Cotlar, C.Sadosky through Krein's moment theory. Their approach (which can
be referred as Generalized Bochner Theorem) provides both integral representation and
extension of forms and kernels invariant under the shift operator. Being applied to a
special bilinear form built with the help of the Hilbert transform and two measures,
this approach gives a necessary and sufficient condition for the Hilbert transform to be
bounded between $L^2$-spaces with respect to these measures (see
[CS1]). The approach of  M.Cotlar, C.Sadosky is very intersting because it provides 
a direct link between the lifting theory of Sz.-Nagy and Foias and thus the scattering
theory) and the continuity of the Hilbert transform in weighted spaces (see [S]).

In the case of equal measures Cotlar-Sadosky theorem is a direct analog of
Helson-Szeg\"{o} characterization of $A_2$ weights.

On the other hand Cotlar-Sadosky theory leaves many questions. It does not provide any
analog of the Hunt-Muckenhoupt-Wheeden characterization of $A_2$ weights. And it seems
impossible to provide in their frame the treatment of general Calder\'{o}n-Zygmund
operator between two weighted spaces.

There is an extensive literature consisting of separately necessary
and sufficient conditions in the spirit of the Hunt-Muckenhoupt-Wheeden characterization
of $A_2$ weights. In most cases authors assume some sort of
$A_{\infty}$ for one (or both) of the weights. Another kind of results is represented very
well by the following theorem of Dechao Zheng in [Z]: $H$ is bounded from $L^2(w^{-1})$ to
$L^2(v)$ if for a positive number $\eta$
$$
\sup_{z \in {\mathbb D}} v^{1 +\eta}(z) w^{1 +\eta}(z) < \infty \, .
$$
Sarason asked a question essentially equivalent to the question whether one can get rid of
$\eta$. The counterexample of Nazarov [N] shows that
$\eta$ cannot be made
$0$. On the other hand, in [TVZ] it is shown that one can weaken the above assumption of
Dechao Zheng if one uses certain Orlic norms instead of $L^{1 + \eta}$.

Here we add another  list of sufficient conditions to the existing collection of such
conditions. We do not know whether our list represent necessary conditions. Most probably it
does not. However, our conditions are completely different of those found before.
We found them by the means of the same Bellman functions we have used above to solve the
two weights problem for the Haar multiplies (which we consider as a ``discrete analog" of two
weights problem for $H$).

Let us introduce further notations. Let $P_z(t) = \frac{1-|z|^2}{|1 - \bar{z}t|}$ be the
Poisson kernel with pole at $z \in \mathbb{D}$. Let $G(z,\zeta) = \log \left| \frac{1 -
\bar{z}\zeta}{z - \zeta}\right|$ be  Green's function with pole at $z$. If $\mu$ is a measure
on $\mathbb{D}$ then Green's potential of this measure is given by $G(\mu)(z)
\stackrel{def}{=} \int_{\mathbb{D}}G(z,\zeta) \,d\mu(\zeta)$. It is positive if $\mu$ is
positive and $\Delta(G(\mu)) = \mu$ in the sense of distributions. We always abbreviate
$G(fdxdy)(z)$ to $G(f)(z)$.

If $u(z)$ is a harmonic function then $u^{'}(z)$ denotes the holomorphic function
$\frac{\partial u}{\partial z}$.

\bigskip

{\bf Lemma 8.1}. {\it Let} $f = (f_1,...,f_k)$ {\it be a} $k$-{\it tuple of real valued}
$L^1(\mathbb{T})$-{\it functions,}
$f(z)$  {\it be a corresponding} $k$-{\it tuple of harmonic functions,} $f^{'}(z)$ {\it be
a corresponding} $k$-{\it tuple of analytic functions.  Let} $B$  {\it be a function of} $k$
{\it real variables, and let} $d^2 B$ {\it denote its Hessian. Put} $b(z) = B(f(z))$. {\it
Then} $\Delta b(z) = 4(d^2 B(f(z)) f^{'}(z), f^{'}(z))$, {\it where} $(\cdot,\cdot)$ {\it
means the scalar product in}
${\mathbb C}^k$.

\bigskip

{\bf Proof}. This is the direct computation of $\frac{\partial^2 B}{\partial \bar{z}\partial
z}$ using the harmonicity of $f(z)$. See [NT].

\bigskip

Our last notation concerns an operator with positive kernel, which is related to $H$ the same
way as $T_{0}$ was related to $T_{\e}$. We first write the bilinear form of this $H_{0}$. Let
$f,g$ be from $L^2(\mathbb{T})$, and let $(\cdot,\cdot)$ {\it
means now the scalar product in}
$L^2(\mathbb{T})$. Then we fefine $H_{0}$ as follows
$$
(H_{0}f,g) \stackrel{def}{=} \int_{\mathbb{D}} f(z)g(z)\left|\frac{w^{'}(z)}{w(z)}\right|
\left|\frac{v^{'}(z)}{v(z)}\right|(1 - |z|) \, dxdy \,.
$$

One can easily see that $H_{0}$ is the operator with positive kernel ($s,t \in \mathbb{T}$)
$$
h(s,t) = \int_{\mathbb{D}} P_{z}(t) P_{z}(s)\left|\frac{w^{'}(z)}{w(z)}\right|
\left|\frac{v^{'}(z)}{v(z)}\right|(1 - |z|) \, dxdy \,.
$$

Now we are ready to formulate and to prove our result about $H$.

\bigskip

{\bf Theorem 8.2}. {\it The following conditions together are sufficient for the boudedness
of} $H$ {\it from} $L^2(w^{-1})$ {\it to} $L^2(v)$:

\begin{eqnarray}
&&\sup_{z \in \mathbb{D}} w(z)v(z) < \infty \, ;\\[.2cm]
&& G(|v^{'}|^2 w)(z) \leq C v(z) \, ;\\[.2cm]
&& G(|w^{'}|^2 v)(z) \leq C w(z) \, ;\\[.2cm]
&& \|H_{0}\|_{L^2(w^{-1}) \to L^2(v)} < \infty 
\end{eqnarray} 

\bigskip

{\bf Remark}. Thus, the problem for singular integral operator is reduced to the problem 
for the operator with positive kernel. In the discrete case considered in Sections 1-7
the corresponding kernel
$k$ did not  satify the regularity conditions usually imposed on the kernels to apply Sawyer's
theory (see [SW]). However, we found  (see Theorem 0.3) necessary and sufficient conditions
for the operator $T_{0}$ with kernel $k$ to be bounded. Now we can do the same. The kernel $h$
of $H_{0}$ does not satisfy in general any regularity conditions unforunately. Still we can
immitate Theorem 0.3 and prove

\bigskip

{\bf Theorem 8.3}. {\it Put} $K(z) = G(|w^{'}v^{'}|)(z)$, {\it and put} $M(z)=G(Kw
\frac{|w^{'}v^{'}|}{wv}),
\, N(z)=G(Kv\frac{|w^{'}v^{'}|}{wv})$.  {\it Then the following conditions together are
sufficient for the boudedness of}
$H$ {\it from} $L^2(w^{-1})$ {\it to} $L^2(v)$:

\begin{eqnarray}
M(z) \leq Cw(z); \\[.2cm]
N(z) \leq Cv(z).
\end{eqnarray}

\bigskip

We leave Theorem 8.3 to the reader. Now we prove Theorem 8.2.

{\bf Proof of Theorem 8.2}. Let $f,g$ be continuous real valued functions on $\mathbb{T}$.
First of all we can assume $g(0)=0$. We also may assume that $w,v$ are continuous because the
estimates we are going to get  are independent from these assumptions.
Clearly,
$$
\left| \int_{\mathbb {T}} Hf \cdot g dm \right| \leq
C\int_{\mathbb{D}}|f^{'}(z)||g^{'}(z)|(1-|z|) \, dxdy \,.
$$

Notice that, immitating the sums from Section 2, we can write

\begin{eqnarray*}
&& \hspace*{-1.1in}\int_{\mathbb{D}}|f^{'}(z)||g^{'}(z)|(1-|z|) \, dxdy  \\[.2cm]
&& \hspace*{-1.1in} \leq 
\int_{\mathbb{D}}|f(z)||g(z)|\left|\frac{f^{'}(z)}{f(z)}-\frac{w^{'}(z)}{w(z)}\right|
\left|\frac{g^{'}(z)}{g(z)}-\frac{v^{'}(z)|}{v(z)}\right|(1-|z|)
\, dxdy \\[.2cm]
&& \hspace*{-1.1in} + 
\int_{\mathbb{D}}|f(z)||g(z)|\left|\frac{w^{'}(z)|}{w(z)}\right|
\left|\frac{g^{'}(z)|}{g(z)}-\frac{v^{'}(z)}{v(z)}\right|(1-|z|)
\, dxdy \\[.2cm]
&&  \hspace*{-1.1in} + 
\int_{\mathbb{D}}|f(z)||g(z)|\left|\frac{v^{'}(z)}{v(z)}\right|
\left|\frac{f^{'}(z)}{f(z)}-\frac{w^{'}(z)}{w(z)}\right|(1-|z|)
\, dxdy \\[.2cm]
&&  \hspace*{-1.1in} + 
\int_{\mathbb{D}}|f(z)||g(z)|\left|\frac{w^{'}(z)}{w(z)}\right|
\left|\frac{v^{'}(z)}{v(z)}\right|(1-|z|)
\, dxdy \\[.2cm]
& = & \Sigma_1 + \Sigma_2 + \Sigma_3 + \Sigma_4 
\end{eqnarray*}

Now let us think that $f=\phi w^{1/2}, \, g =\psi v^{1/2}$. We only need to prove that
$\Sigma_i, \, i=1,2,3$ are estimated by $C\|\phi\|_2^2 + \|\psi\|_2^2$. In fact,replacing
$f$ by $tf$ and $g$ by $t^{-1}g$ allows then to estimate these sums by $C\|\phi\|_2
\|\psi\|_2$.

Let us first do it for $\Sigma_1$. Consider $B(X,x,w,Y,y,v) \df X-\frac{x^2}{w} +
Y-\frac{y^2}{v}$. Put $b(z) = B(\phi^2(z), \phi w^{1/2}(z), w(z), \psi^2(z), \psi
v^{1/2}(z), v(z))$. Using Lemma 8.1 we can compute its Laplacian and see that it is
superharmonic and moreover that ($f(z) =\phi w^{1/2}(z), \, g(z)=\psi
v^{1/2}(z)$)

$$
-\Delta(b(z))  \geq 
c\frac{f(z)^2}{w(z)}\left|\frac{f^{'}(z)}{f(z)}-\frac{w^{'}(z)}{w(z)}
\right|^2 +
c\frac{g(z)^2}{v(z)}\left|\frac{g^{'}(z)}{g(z)}-\frac{v^{'}(z)}{v(z)}\right|^2
$$
Thus,
$$
-\Delta(b(z))   \geq   
c\frac{|f(z)g(z)|}{\sqrt{(w(z)v(z)}}\left|\frac{f^{'}(z)}{f(z)}-\frac{w^{'}(z)}{w(z)}
\right|\left|\frac{g^{'}(z)}{g(z)}-\frac{v^{'}(z)}{v(z)}\right|
$$

And so we obtain using (8.1):

\begin{equation}
-\Delta(b(z)) \geq c|f(z)g(z)|\left|\frac{f^{'}(z)}{f(z)}-\frac{w^{'}(z)}{w(z)}
\right|\left|\frac{g^{'}(z)|}{g(z)}-\frac{v^{'}(z)}{v(z)}\right|
\end{equation}

Notice also that function $b$ vanishes on the circle $\mathbb{T}$. In fact,
$\phi^2(z)=\frac{(\phi w^{1/2}(z))^2 }{w(z)}, \psi^2(z)= \frac{(\psi v^{1/2}(z))^2}{v(z)}$ on
$\mathbb{T}$. Notice also that
$b(0) \leq \|\phi\|_2^2 + \|\psi\|_2^2$.

Applying Green's formula, we get

\begin{equation}
-\int_{\mathbb{D}} \Delta(b(z))\log\frac{1}{|z|} \,dxdy = b(o) - \int_{\mathbb{T}} b \,dm \leq
\|\phi\|_2^2 + \|\psi\|_2^2
\end{equation}

Combining (8.7) and (8.8), we obtain that
$$
\Sigma_1 \leq C\|\phi\|_2^2 + \|\psi\|_2^2 \,.
$$

Estimate of $\Sigma_2, \Sigma_3$ are similar to each other. 

Let us estimate$\Sigma_2$. To do that we need the following $p(z)$ (remind that $f=\phi
w^{1/2}$) :
$$
p(z) \df \phi^2(z) - \frac{f(z)^2}{w(z) + G(|w^{'}|^2 v)(z)} + \psi^2(z) -
\frac{g(z)^2}{v(z)} \,.
$$
Notice how small is the difference with $b$. But this small difference allows the following
estimate of the Laplacian (we use here Lemma 8.1 again and also formula (6.10)):

$$
-\Delta(p(z)) \geq  c\frac{f(z)^2}{(w(z) + G(|w^{'}|^2 v)(z))^2}\Delta(w(z) + G(|w^{'}|^2
v)(z)) + c\frac{g(z)^2}{v(z)}\left|\frac{g^{'}(z)}{g(z)}-\frac{v^{'}(z)}{v(z)}\right|^2
$$ 
Thus,
$$
-\Delta(p(z)) \geq
c\frac{|w^{'}(z)|^2}{w(z)^2} f(z)^2 v(z) +
c\frac{g(z)^2}{v(z)}\left|\frac{g^{'}(z)}{g(z)}-\frac{v^{'}(z)}{v(z)}\right|^2
$$

We used here our second assumption from Theorem 8.2: $G(|w^{'}|^2 v)(z) \leq Cw(z)$.
Thus,

\begin{equation*}
-\Delta(p(z)) \geq  
c\frac{|w^{'}(z)|}{w(z)} |f(z) 
g(z)|\left|\frac{g^{'}(z)}{g(z)}-\frac{v^{'}(z)}{v(z)}\right|
\end{equation*}

Again $p$ vanishes on $\mathbb{T}$ and again obviously $p(0) \leq \|\phi\|_2^2 +
\|\psi\|_2^2$.

Thus, using Green's formula as before, we get

$$
\int_{\mathbb{D}}\frac{|w^{'}(z)|}{w(z)} |f(z) 
g(z)|\left|\frac{g^{'}(z)}{g(z)}-\frac{v^{'}(z)}{v(z)}\right| \log\frac{1}{|z|} \, dxdy
\leq
c^{-1}(\|\phi\|_2^2 +\|\psi\|_2^2),
$$
which is the estimate of $\Sigma_2$. Similarly we treat $\Sigma_3$ using the third assumption
of Theorem 8.2: $G(|v^{'}|^2 w)(z) \leq Cv(z)$.
Theorem 8.2 is completely proved.

\setcounter{equation}{0}

\section{9.Two weight norm inequalities for S-functions. Necessary and sufficient
conditions.}

Square functions play an important role in the theory of Singular integral operators. Often
the estimate of Singular integral operators goes through the estimate of a certain
$S$-function. The multitude of examples can be found in [St]. In one of the most recent
examples in [V] this approach was applied to characterize the matrix ${\bf A}_p$ weights.

Let us consider
$$
S(f)(x) \df \sum_{I: x \in I}|\langle f\rangle_{I_{-}} - \langle f\rangle_{I_{+}}|^2 \,.
$$

So when $f \rightarrow S(f)$ is bounded from $L^2(w^{-1})$ to $L^2(v)$? The following theorem
gives the answer. It says basically that this is so if and only if there is a uniform bound
on test functions $f= w \chi_J$.

\bigskip

{\bf Theorem 9.1}.{\it $\|S(f)\|_{L^2(v)} \leq C\|f\|_{L^2(w^{-1})}$ if and only if (1.1) and
(1.2) hold.}

\bigskip

{\bf Proof}. Necessity is simple. To prove the sufficiency we again use the Bellman function
approach. Let us consider 
$$
B(X, x, w, M) = CX - \frac{x^2}{w} - \frac{x^2}{w + M} \,.
$$
We are going to compose it with the following (super)martingales
$$
X_I=\langle\phi^2\rangle_{I}, \, x_I = \langle \phi w^{1/2}\rangle_{I}, \,  w_I =\langle w
\rangle_I
$$
$$
M_I = \frac{1}{|I|}\sum_{\ell \subseteq I}|\langle w\rangle_{\ell_{-}} - \langle
w\rangle_{\ell_{+}}|^2 \langle v\rangle_{\ell}|\ell| \,.
$$
By (1.2) $M_I \leq C w_{I}$. Then we have the supermartingale
$$
b(I) = B(X_I, x_I, w_I, M_I) \,.
$$
Its ``discrete Laplacian" $\Delta (b)(I) \df \frac{b(I_{-})+b(I_{+})} {2} - b(I)$ can be
estimated using Lemma 9.2 bellow (we use the notations
$f = \phi w^{1/2}$):

$$
-\Delta (b)(I) =  b(I) - \frac{b(I_{-})+b(I_{+})} {2} \geq 
$$
$$
c\frac{\langle
f\rangle_I^2}{\langle w\rangle_I}
\inf_{c_1,c_2 \in [1/2,2]}\left|c_1\frac{\langle f\rangle_{I_{-}} -
\langle f\rangle_{I_{+}}}{\langle f\rangle_{I}} - 
c_2
\frac{\langle w\rangle_{I_{-}} - \langle w\rangle_{I_{+}}}{\langle w\rangle_{I}}\right| +
$$
$$
c\frac{\langle
f\rangle_I^2}{\langle w\rangle_I^2} |\langle w\rangle_{I_{-}} - \langle w\rangle_{I_{+}}|^2
v_I
$$
We want to continue the estimate of the negative Laplacian from below. To this end let us
consider two cases.

{\bf First case}: $\left|\frac{\langle w\rangle_{I_{-}} - \langle w\rangle_{I_{+}}}{\langle
w\rangle_{I}}\right| \leq \frac{1}{10}\left|\frac{\langle f\rangle_{I_{-}} - \langle
f\rangle_{I_{+}}}{\langle f\rangle_{I}}\right|$. In this case we use the first term in the
estimate of our negative Laplacian to see that
$$ 
-\Delta (b)(I)  \geq 0.09c \cdot \frac{\langle
f\rangle_I^2}{\langle w\rangle_I}\left|\frac{\langle f\rangle_{I_{-}} - \langle
f\rangle_{I_{+}}}{\langle f\rangle_{I}}\right|^2  
$$
And so,

\begin{equation}
-\Delta (b)(I) \geq c^{'} \frac{|\langle f\rangle_{I_{-}} - \langle
f\rangle_{I_{-}}|^2}{\langle w\rangle_I} \geq c^{''}|\langle f\rangle_{I_{-}} - \langle
f\rangle_{I_{+}}|^2 \langle v\rangle_I \,.
\end{equation}
In the last inequality we use (1.1).

{\bf Second case}: $\left|\frac{\langle w\rangle_{I_{-}} - \langle w\rangle_{I_{+}}}{\langle
w\rangle_{I}}\right| \geq \frac{1}{10}\left|\frac{\langle f\rangle_{I_{-}} - \langle
f\rangle_{I_{+}}}{\langle f\rangle_{I}}\right|$. In this case we use the second term in the
estimate of the negative Laplacian. Taking into account the ineqality above we get

\begin{equation}
-\Delta (b)(I)  \geq 0.01 \cdot|\langle f\rangle_{I_{-}} - \langle
f\rangle_{I_{+}}|^2 \langle v\rangle_I \,.
\end{equation}

So we get the same estimate as in the first case. Now we can just apply Green's formula for
our discrete Laplacian, which amounts to just applying (9.1),(9.2) to $I$, then to
$I_{\pm}$, then to ``sons" of $I_{\pm}$,  et cetera... to obtain

$$
\frac{1}{|J|}\sum_{I \subseteq J} |\langle f\rangle_{I_{-}} - \langle
f\rangle_{I_{+}}|^2 \langle v\rangle_I |I| \leq C B(X_J, x_J, w_J, M_J) \leq C \langle
\phi^2\rangle_J
\,.
$$
But $\langle \phi^2\rangle_J = \|f\|^2_{L^2(w^{-1})}$, and Theorem 9.1 is completely proved.

\bigskip

We are left to prove the following lemma. Let us denote $a =(X,x,w)$, and let $ P(a) \df X -
\frac{x^2}{w}, \, Q(a,M) \df X - \frac{x^2}{w+M}$.

\bigskip

{\bf Lemma 9.2}. {\it Let $a = \frac{a_{-}+a_{+}}{2}, \, M \geq \frac{M_{-}+M_{+}}{2}$.
Then}
$$
-d^2 P = 2\frac{x^2}{w}\left|\frac{dx}{x} - \frac{dw}{w}\right|^2 \, ;
$$
$$
Q(a,M) - \frac{Q(a_{-})+Q(a_{+})}{2} \geq c\frac{x^2}{w^2}(M - \frac{M_{-}+M_{+}}{2})\, ;
$$
$$
P(a) - \frac{P(a_{-})+P(a_{+})}{2} \geq c \frac{x^2}{w} \inf_{c_1,c_2
\in [1/2,2]}\left|c_1\frac{x_{-}-x_{+}}{x} -
c_2\frac{w_{-}-w{+}}{w}\right|^2 \, .
$$

\bigskip

{\bf Proof}. The first equality is obtained by direct computation. It shows that $B$ is
concave. Let $a(t)$ be the linear function on $[-1,1]$ assuming values $a_{-},a_{+}$ at
endpoints (and thus $a_0 = a$). Let $M(t)$ be the  piecewise linear function on $[-1,1]$
assuming values
$M_{-},M,M_{+}$  at $-1,0,1$ correspondingly. Then $M(t)$ is concave, and $M^{''}(t) = 
(\frac{M_{-}+M_{+}}{2} - M)\delta_0 $, where $\delta_0$ is a Dirac measure at $0$. Notice
that $Q$ is a composition of $P$ and a linear function. So $Q$ is concave. Thus, by (6.10)
$-d^2 Q \geq -\frac{\partial Q}{\partial M} d^2 M$. In particular, if $q(t) \df Q(a(t),M(t))$
then measure $q^{''}(t)$ satisfies
$$
-q^{''}(t) \geq c\frac{x(t)^2}{w(t)^2}(M - \frac{M_{-}+M_{+}}{2})\delta_0 \, .
$$
Thus
$$
q(0)- \frac{q(-1) + q(1)}{2} = -\int_{-1}^1 (1-|t|)q^{''}(t) \geq c\frac{x^2}{w^2}(M -
\frac{M_{-}+M_{+}}{2})\, ,
$$
and the second inequality is proved.

To prove the third inequality of the lemma let us introduce $p(t) \df P(a(t))$. Using the
calculation for $-d^2P$ we obtain that then measure $p^{''}(t)$ satisfies
$$
-p^{''}(t)  = 2\frac{x(t)^2}{w(t)}\left|\frac{x_{-}-x_{+}}{x(t)} -
\frac{w_{-}-w_{+}}{w(t)}\right|^2 \,.
$$

On $[-1/2,1/2]$ we have that $\frac{x(t)}{x} \in[1/2,2]$ and $\frac{w(t)}{w} \in[1/2,2]$.
Thus, on $[-1/2,1/2]$, we have
$$
-p^{''}(t)  \geq  c\frac{x^2}{w}\inf_{c_1,c_2 \in [1/2,2]}\left|c_1\frac{x_{-}-x_{+}}{x} -
c_2\frac{w_{-}-w_{+}}{w}\right|^2 \,.
$$
We finish the proof by combining this inequality with the following one

$$
p(0)- \frac{p(-1) + p(1)}{2} = -\int_{-1/2}^{1/2} (1-|t|)p^{''}(t) \,.
$$

\section{Concluding remarks}

 1) Seems like the case $p=2$ is a true miracle, because in this
case we were
able to give a finite list of simple conditions which are necessary and
sufficient for
two weight boundedness of our family $T_{\pm}$ of Calder\'{o}n-Zygmund
operators. As for $p \neq 2$ case, there are strong indications that the similar list of
conditions (which can be actually copied from the ``$p=2$'' case) will not be equivalent
to two weights estimate. 

Let us notice that there exists an approach through Cotlar-Sadosky theory to
two weight estimate for the Hilbert transform even for $p \neq 2$ (see [CS2]).

\bigskip

2) One wonders whether the 5 conditions in our list in Section 1 are independent. This
is most probably so, but
the proof should be quite involved. At least it follows from [N] that (1.1) does not imply
neither (1.2) nor (1.3).

\end{document}